\documentclass[12pt]{amsart}
\usepackage{amssymb}
\usepackage[active]{srcltx}
\usepackage{epsfig}
\usepackage{amscd}
\usepackage{a4wide}
\usepackage{amsthm}
\usepackage{enumitem}
\usepackage{tikz}
\usepackage{tikz-cd}
\usetikzlibrary{matrix}
\usepackage{pgfplots}
\usepackage{hyperref}
\usepackage{derivative}
\usepackage{algorithm}
\usepackage[noend]{algpseudocode}
\usepackage{array}
\usepackage{multirow}

\usepackage{caption}
\usepackage{subcaption}
\usepackage{adjustbox}

\newcommand{\wcf}[2]{\ensuremath[#2]_{#1}}
\newcommand{\swcf}[1]{\ensuremath[#1]}

\newtheorem{theorem}{Theorem}[section]
\newtheorem{lemma}[theorem]{Lemma}
\newtheorem{proposition}[theorem]{Proposition}
\newtheorem{corollary}[theorem]{Corollary}

\theoremstyle{definition}
\newtheorem{definition}[theorem]{Definition}

\newtheorem{remark}[theorem]{Remark}

\numberwithin{equation}{section}

\newcommand{\bc}{{\mathbb C}}
\newcommand{\bz}{{\mathbb Z}}
\newcommand{\br}{{\mathbb R}}
\newcommand{\abs}[1]{\left\lvert{#1}\right\rvert}
\newcommand{\wcfM}[1]{M_{#1}}
\let\Re\relax
\DeclareMathOperator{\Re}{Re}
\let\Im\relax
\DeclareMathOperator{\Im}{Im}

\def\ns {\mathbf{n}}
\def\nsL {\mathbf{n}_L}
\def\nsR {\mathbf{n}_R}

\newcommand{\im}{{\sqrt{-1}}}
\newcommand{\inprod}[2]{\langle {#1}, {#2} \rangle}

\definecolor{geoblue}{rgb}{0.08235294117647059,0.396078431372549,0.7529411764705882}
\definecolor{geored}{rgb}{0.8274509803921568,0.1843137254901961,0.1843137254901961}

\newtheorem{innercustomthm}{Theorem}
\newenvironment{customthm}[1]
  {\renewcommand\theinnercustomthm{#1}\innercustomthm}
  {\endinnercustomthm}

\DeclareMathOperator{\mindeg}{mindeg}

\begin{document}
\title [Non-freeness of parabolic two-generator groups] {Non-freeness of parabolic two-generator groups}
\author{Philip Choi, Kyeonghee Jo, Hyuk Kim, and Junho Lee}
\address{Department of Mathematical Sciences, Seoul National University, Seoul, 08826, Korea}
\email{philip94@snu.ac.kr}
\address{Division of Liberal Arts and sciences,  Mokpo National Maritime University, Mokpo, Chonnam, 530-729, Korea  }
\email{khjo@mmu.ac.kr}
\address{Department of Mathematical Sciences, Seoul National University, Seoul, 08826, Korea}
\email{hyukkim@snu.ac.kr}
\address{Department of Mathematical Sciences, Seoul National University, Seoul, 08826, Korea}
\email{abab9579@snu.ac.kr}
\thanks {The second author acknowledges support from the National Research Foundation of Korea (NRF) grant funded by the Korean government (MSIT) (NRF-2021R1F1A1049444) and the Open KIAS Center at the Korea Institute for Advanced Study. The third author acknowledges support from the National Research Foundation of Korea (NRF) grant funded by the Korean government (MSIT) (NRF-2018R1A2B6005691).}

\begin{abstract}
A  complex number $\lambda$ is said to be non-free if  the subgroup of $SL(2,\bc)$ generated by 
$$X=\begin{pmatrix}
1& 1\\
0 & 1 
 \end{pmatrix} 
\,\, \text{and}\,\,\,Y_{\lambda}=\begin{pmatrix}
1& 0\\
\lambda & 1 
 \end{pmatrix}$$ is not a free group
of rank 2.
In this case the number $\lambda$ is called a relation number, and it has been a long standing problem to determine the relation numbers. In this paper, we characterize the relation numbers by establishing  the equivalence between $\lambda$ being a relation number and $u:=\sqrt{- \lambda}$ being a root of a `generalized Chebyshev polynomial'.
The generalized Chebyshev polynomials of degree $k$ are given by a sequence of $k$ integers $(n_1, n_2,\cdots, n_k)$ using the usual recursive formula, and thereby can be studied systematically using continuants and continued fractions. Such formulation, then, enables us to prove that, 
the question whether a given number $\lambda$ is a relation number of $u$-degree $k$ can be answered by checking only finitely many generalized Chebyshev polynomials.
Based on these theorems, we design an algorithm 
deciding any given number is a relation number with minimal degree $k$. With its computer implementation we provide a few sample examples, with a particular emphasis on the well known conjecture that every rational number in the interval $(-4, 4)$ is a relation number.
\end{abstract}

\maketitle


\section{Introduction}
For each  complex number $\lambda$, we consider a subgroup of $SL(2,\bc)$ generated by 
$$X=\begin{pmatrix}
1& 1\\
0 & 1 
 \end{pmatrix} 
\,\, \text{and}\,\,\,Y_{\lambda}=\begin{pmatrix}
1& 0\\
\lambda & 1 
 \end{pmatrix}$$
and denote it by $G_{\lambda}$. 
Note that for any two parabolic elements $A, B$ in $SL(2,\bc)$ such that $AB\neq BA$, there exists $U\in SL(2,\bc)$ such that $UAU^{-1}=X$, $UBU^{-1}=Y_{\lambda}$ with $\lambda=2-tr(AB^{-1})$, see  \cite{Riley4}. Therefore any two non-abelian groups $G_1$ and $G_2$ generated by two parabolic elements $X_i$ and $Y_i$ (with $i=1,2$) are conjugate in $SL(2,\bc)$, and hence isomorphic, if $tr(X_1{Y_1}^{-1})$=$tr(X_2{Y_2}^{-1})$. The problem of whether $G_{\lambda}$ is free or not has been the subject of extensive studies by many authors.  
We say that $\lambda$ is {\it free} when $G_{\lambda}$ is isomorphic to a free
group of rank $2$. Otherwise, we say that $\lambda$ is {\it non-free} or a {\it relation number}. (It is unfortunate that our $\lambda$ is twice the “$\lambda$” which is used by other authors in most previous works.)

By definition, $\lambda$ is non-free  if and only if there exists a nontrivial word of the form 
$$
w(x,y)=x^{n_1}y^{n_2}x^{n_3}y^{n_4}\cdots x^{n_{m}}, n_i\in\mathbb Z
$$
in the free group $F=\langle x,y\rangle$ which evaluates to $Id\in SL(2,\bc)$ for $w(X,Y_{\lambda})$. In other words, a non-free number $\lambda$ must satisfy the four polynomial equations
$$w(X,Y_{\lambda})_{11}=1, w(X,Y_{\lambda})_{12}=0, w(X,Y_{\lambda})_{21}=0, w(X,Y_{\lambda})_{22}=1.$$
However, it turns out that  $\lambda$ is a relation number if and only if  $\sqrt{-\lambda}$ is a root of a generalized Chebyshev polynomial (see Theorem A below). A polynomial characterization of the relation number has been considered previously by many authors \cite{LU, Riley4, brenner1975non, Bamberg, Ilia}, and Chebyshev formulation appeared recently for the representation of the 2-bridge knot group case \cite{Jo-Kim2}. It is interesting to notice that the coefficients of Bamberg polynomial has ‘cyclic symmetry’ while those of Chebyshev has ‘linear symmetry’. In this paper we further show that all the components $w(X,Y)_{ij}$ can be written in the Chebyshev polynomials explicitly (see Theorem \ref{W-eq}), which should be useful for a ‘membership problem’ \cite{chorna2017two, Sanov, Brenner}. The Chebyshev setting brought us lots of advantages in the analysis for the proofs of Theorem B and C below. 

In 1955, Brenner \cite{Brenner} demonstrated that $G_{\lambda}$ is free for $|\lambda|\geq 4$ following \cite{Sanov}, then  in 1958, Chang-Jennings-Ree \cite{CJR} extended this by proving that $G_{\lambda}$ is free for $\lambda$ outside three open discs of radius $2$ centered at $(-2,0)$, $(0, 0)$, and $(2, 0)$ respectively. In 1969, Lyndon and Ullman \cite{LU} further contributed by establishing that $G_{\lambda}$ is free if either of the following two conditions hold:
\begin{itemize}
\item  $\lambda$ is not within the interior of the convex hull of the set consisting of the circle $|z|=2$ together with the two points $z=\pm 4$.
\item  $|\lambda\pm i|\geq 1$ and $|\lambda\pm2|\geq 2$.
\end{itemize}
After that, Ignatov proved  that $G_{\lambda}$ is free if 
any one of the following conditions holds (see \cite{ja1976free} and  \cite{Ignatov}):
\begin{itemize}
\item  $|\lambda-2|>1$ and $2\leq |Re (\lambda) |<\frac{5}{2}$.
\item 
$|\lambda|>2 $ and $|Im (\lambda) |>1$.
\end{itemize}
 
Now it is well known that $\lambda$ is free (and discrete) when it lies in the closure of the Riley slice in $\mathbb C$ \cite{Keen-Series} and the set of relation numbers is dense in the exterior of the Riley slice ‘in a strong sense’ \cite{martin2021nondiscrete}. (See also \cite{akiyoshi2021classification} for a recent summary for the discrete case.)

There have also been numerous results about the rational relation numbers. 
It is known earlier that the set of relation numbers are densely distributed inside  the unit circle \cite{Ree} and also dense on the interval $(-4, 4)$ \cite{CJR}. Then one might ask about the non-freeness of rational numbers in that interval, and Lyndon and Ullman actually verified that a family of rational numbers and some others are non-free. Since then the list grows gradually by the works of many authors, and 
 there are currently no known rational numbers within the interval $(-4, 4)$ being identified as free. (see Gilman's paper \cite{Gilman2} for a summary up to that time.) Hence it has been conjectured that all rational numbers within the interval $(-4, 4)$ are non-free (see \cite{Farbman, brenner1975non, KK, Jang-Kim, Ilia, chorna2017two}, and references therein). 
  Notably in 1995, Farbman proved that all rational numbers in  $(-4, 4)$ whose numerator is at most $16$ are non-free \cite{Farbman}, and then in 2022, Kim and Koberda expanded this result to include rational numbers with numerators up to $27$ with the possible exception of $24$ \cite{KK}. In this paper we will consider the rational numbers with fixed denominators up to $22$ as an application of our result. 


It is known to be difficult to find accumulation points of non-free rational numbers. In fact,  there is no known sequence of non-free rational numbers converging to $3$ or $4$.  In 1993, Beardon found a sequence of rational numbers $\{\mu_k\}$ for any non-square integer $n$, such that each $\mu_k$ is a convergent of $\sqrt{n}$ satisfying Pell's equation, and  $\{\frac{1}{\mu_k^2}\}$ is a sequence of non-free rational numbers converging to $\frac{1}{n}$ \cite{beardon1993pell}. Subsequently, in 1996,  Tan-Tan found more sequences of non-free rational numbers  which converge to numbers of the form $\frac{n}{n\pm 1}$ using quadratic Diophantine equations \cite{tan1996quadratic}. In 2021, Smilga discovered sequences  of non-free rational numbers converging to $(\frac{\sqrt{5}+1}{2})^{\pm 2}$ and many other interesting sequences of non-free rational numbers \cite{Ilia}. These sequences are of a very special type, and we can show that such limit points, including $(\frac{\sqrt{5}+1}{2})^{\pm 2}$, are also relation numbers  (see Remark \ref{Smilga}).

In Section \ref{ChebyContinuant}, we introduce generalized Chebyshev polynomials, and discuss various properties of these polynomials.
In Section \ref{section:main-A}, we prove that a relation number $\lambda$ can be characterized by the property that $\sqrt{- \lambda}$ is a root of a certain class of  `generalized Chebyshev polynomial' as follows.
\begin{customthm}{A}[Theorem \ref{rel number}, Corollary \ref{3rd-4th}, Theorem \ref{factor-cont}] \label{main2}
Let $\lambda$ be a nonzero complex number.
\begin{enumerate}
\item The following are all equivalent:
\begin{enumerate}
\item[\rm (i)]
$\lambda$ is a relation number.
\item[\rm (ii)]
$\sqrt{-\lambda}$ is a root of an $\bf n$-Chebyshev polynomial $s^{\bf n}_{k+1}(t) \in \mathbb Z[t]$ for some nonzero integer sequence ${\bf n}=(n_1,n_2,\cdots,n_k)$.
\item[\rm (iii)]
$\sqrt{-\lambda}$ is a root of an $\bf n$-Chebyshev polynomial $s^{\bf n}_{2k}(t) \in \mathbb Z[t]$ for some nonzero integer sequence ${\bf n}=(n_1,n_2,\cdots,n_{2k-1})$.
\item[\rm (iv)] $\sqrt{-\lambda}$ is a root of a polynomial $T^{\bf n}_{k}(t)\in \mathbb Z[t]$
 for some nonzero integer sequence ${\bf n}=(n_1,n_2,\cdots,n_k)$.
\item[\rm (v)] $\sqrt{-\lambda}$ is a root of a polynomial $\tilde{T}^{\bf n}_{k}(t)\in \mathbb Z[t]$  for some nonzero integer sequence ${\bf n}=(n_1,n_2,\cdots,n_k)$.
\item[\rm (vi)] $\sqrt{-\lambda}$ is a root of a polynomial $\hat{T}^{\bf n}_{k}(t)\in \mathbb Z[t]$
 for some nonzero integer sequence ${\bf n}=(n_1,n_2,\cdots,n_k)$.
\end{enumerate}
\item $\lambda$ is a relation number if  $\sqrt{-\lambda}$ is a zero of either an $\bf n$-Chebyshev polynomial $v^{\bf n}_{k}(t)$ or 
$v^{\bf n}_{k}(-t)$ in $\mathbb Z[t]$ for some nonzero integer sequence ${\bf n}=(n_1,n_2,\cdots,n_k)$.
\end{enumerate}
\end{customthm}
Here $s^{\mathbf{n}}_{k+1}(t)$ is a polynomial of degree $k$ whose coefficients are constructed from the sequence ${\bf n}=(n_1,n_2,\cdots,n_k)$ using Chebyshev polynomial recursion generalizing the classical Chebyshev polynomials of the second kind, while  $T^{\mathbf{n}}_{k}(t)$  is a generalized version of the classical Chebyshev polynomial of the first kind with degree $k$. $\tilde{T}^{\mathbf{n}}_{k}(t)$ and $\hat{T}^{\mathbf{n}}_{k}(t)$ are slight variations of $T^{\mathbf{n}}_{k}(t)$. Furthermore, 
$\tilde{T}^{\bf n}_{2k+1}(t)$ coincides with Bamberg's polynomial $B_k$ in \cite{Bamberg} as $(-1)^ktB_k(-\frac{t^2}{2})$, and $\hat{T}^{\bf n}_{k+1}(t)$ with Smilga's polynomial $P_{HR}^{k+1}(t)$ in \cite{Ilia} as either  $tP_{HR}^{k+1}(\mathbf{n};-t^2)$ (when $k$ is odd) or $t^2P_{HR}^{k+1}(\mathbf{n};-t^2)$ (when $k$ is even). Also,  Jang-Kim's polynomial in  \cite{Jang-Kim} is a specialization of the $T^{\bf n}_{2k+1}(t)$, where ${\bf n}=(1,-1,1,-1,\ldots)$, divided by a factor of $t$. (see Remark \ref{Jang-Kim}).
The polynomials $v^{\bf n}_{k}(t)$ and  $v^{\bf n}_{k}(-t)$ in (2), are also generalized Chebyshev polynomials corresponding to the third and the fourth kind,
respectively.  (See Section \ref{ChebyContinuant} and \ref{section:main-A} for the details.)



As easy consequences of Theorem \ref{main2}, we will demonstrate in Section \ref{section:rational}, through a direct factorization of low degree Chebyshev  polynomials $s^{\bf n}_{2k}$ for $k\leq 4$, that the following explicitly given rational numbers all become the relation numbers   for any nonzero integers $r, s, t, v$:
\begin{itemize}
   \item $\displaystyle\ \ \frac{r+s}{rs},\,  \frac{r+s+t}{rst}, \frac{r^2s+tv^2}{r^2stv+rstv^2},\,\, \frac{1}{r}+\frac{1}{rs}+\frac{1}{rst}$
   \item $\displaystyle\ \ \frac{r+w+t}{rst}\quad \text{when}\,\,\,w|rs$
   \item $\displaystyle\ \ \frac{r+t+v+w}{rst}\quad \text{when}\,\,\, vw\,|\,rst,\,\, rv+tv+tw=0$
   \item $\displaystyle\ \ \frac{r+t+v+w}{rst}\quad \text{when}\,\,\, rst=vw, \,\,r^2s+w^2+rst=(r+t)(v+w)$
\end{itemize}

Of course we can find lots of rational relation numbers directly from these lists but these do not cover all the rational numbers in $(-4,4)$ since they are obtained only from $s^{\bf n}_{2m}$ for $m\leq 4$, which are determined by an integer sequence of length $7$, ${\bf n}=(n_1,\cdots,n_7)$. In fact, we can show that as a rational number approaches $\pm 4$, we need higher and higher degree Chebyshev polynomials. (See Proposition \ref{lower bound min deg}.)  

To resolve the rational relation number conjecture, for a given number $\lambda$ we need a way to check that $\sqrt{-\lambda}$ is a root of some Chebyshev  polynomial $s^{\bf n}_{k+1}$, ${\bf n}=(n_1,n_2,\cdots,n_k)$. And this can be done  in principle when such $k$ is minimal by the following `finiteness theorem' not for only rational but for any complex number $\lambda$.
\begin{customthm}{B}\label{main-higher step}
Let $\lambda = -u^2 \in \mathbb{C}$ be a relation number, and let $k$ be the smallest positive integer such that $s^{\bf n}_{k+1}(u) = 0$ for some integer sequence ${\bf n}=(n_1, n_2, \cdots, n_{k})$. Then there exist positive integers $N_{1}, N_{2}, \cdots, N_{k}$ such that $|n'_i| \leq N_{i}$ ($ i=1,\cdots,k$) holds for all integer sequences ${\bf n'}$ with $s^{\bf n'}_{k+1}(u) = 0$. 
\end{customthm}


Note that the above theorem does not hold for a non-minimal  $k$. That is, if there is a smaller positive integer $k'$ than $k$ such that $s^{\bf n}_{k'+1}(u) = 0$  for some integer sequence  ${\bf n}$, then there might be infinitely many sequences ${\bf n}$ such that $s^{\bf n}_{k+1}(u) = 0$. For example, 
while $(1,1)$ and $(-1, -1)$ are the only two integer sequences of length $k = 2$ satisfying $s^{\bf n}_3(1)=0$,
infinitely many sequences of length $5$, including ${\bf n} = (1,1,n,\pm 1,\pm 1)$ for all $n\in \mathbb Z$, satisfy $s^{\bf n}_{6}(1) = 0$.

Hence the problem of deciding a given number is a relation number of minimal length or $u$-degree $k$ can be solvable at least theoretically once we establish an algorithm to find $N_1,\cdots,N_k$. Indeed, we have the following `decision theorem'. 
\begin{customthm}{C}\label{decision algorithm}
  There is an algorithm to decide if a given $\lambda=-u^2\in \mathbb{C}$ is a relation number with a given minimal $u$-degree $k$.
\end{customthm}

If we denote the  set of all relation number with minimal $u$-degree $k$ by $R_k$, then each relation number belongs to $R_k$ for some positive integer $k$, and the set of relation numbers is the disjoint union of $R_k$’s. Theorem \ref{decision algorithm} answers the decision problem of whether a given complex number is contained in $R_k$ or not. 



The proof of Theorem \ref{main-higher step} is given in Section \ref{section:pf main-higher step}, and then the algorithm for decision problem, Theorem \ref{decision algorithm}, follows in the next section. The initial idea for the proof is the observation that, as some of $n_i$'s approach the infinity, the corresponding roots of Chebyshev polynomials approach the roots of some lower degree Chebyshev polynomials.

In Section \ref{section:algorithm}, we prove Theorem \ref{decision algorithm} by establishing an exhaustive algorithm to find a nonzero integer sequence ${\bf n}$ of a given length $k$ satisfying $s^{\bf n}_{k+1}(\sqrt{-\lambda}) = 0$ for a given $\lambda \in \mathbb{C}$ assuming $k$ is the minimal such length.
The idea for the proof is to consider a (generalized) continued fraction associated to $\lambda$ and ${\bf n}$. 
We show that the shortest continued fraction whose value is $\infty$ gives the sequence we want to find. The algorithm is designed to systematically compute the maximum absolute value of continued fractions of a given length $i$, denoted as $M_i$, and it checks whether $M_i$ equals to $\infty$ or not. The values of $M_i$ are calculated in sequential order, which establishes that $M_i$ is achieved by an integer sequence within a bounded region of $\mathbb{Z}^i$ determined using $M_1, M_2, \cdots, M_{i-1}$. These also essentially determine $N_i$ in Theorem \ref{main-higher step}. 

But in the actual computer implementation,  as $i$ increases, the size of the region grows exponentially. And therefore we introduce a greedy algorithm that, in conjunction with the exhaustive algorithm, enables us to find the sequences more quickly, although they might not be of the minimal length. We tried to list the integer sequences for all the rational numbers in the interval $(-4,4)$ with denominators up to $22$.
We determine all such sequences except possibly for $6$ fractions:
\[
\frac{51}{13}, 
\frac{67}{17},
\frac{74}{19}, \frac{75}{19},
\frac{83}{21}, 
\frac{87}{22}.
\]
It takes too much time (more than a week with a desktop computer) for these numbers, and we could not conclude whether these are relation numbers or not with our computer. 

We provide the corresponding integer sequence ${\bf n}$ for each of these rational numbers with denominator $\leq 22$, along with a usable Mathematica code in \url{https://relation-numbers.diagram.site/}. 
We only present the computed results when the denominator is $13$ as a sample in this paper. In addition, a few irrational real and complex relation numbers are also presented using this algorithm.



\section{${\bf n}$-Chebyshev polynomials and Continuants}\label{ChebyContinuant}
In this section, we define an $\bf n$-Chebyshev polynomial  for an arbitrary integer sequence ${\bf n}=(n_1,n_2,\cdots, n_k)$, as a generalization of Chebyshev polynomials, and investigate their  relationship with the continuants. The result of this section will be used later.
\subsection{Chebyshev polynomials}\label{Cheby}
$Chebyshev$ $polynomials$ are defined by a three-term recursion
$$g_{n+1}(t)=tg_n(t)-g_{n-1}(t).$$ 
We denote the Chebyshev polynomials $g_n(t)$ with the initial condition 
$g_0(t)=a, g_1(t)=b$ by $Ch_n^t(a,b)$, which depends on the initial condition linearly.  And the following properties are also obvious:
\begin{equation}\label{linearity of Cheby}
\begin{split}
Ch_n^t(a,b)&=Ch_{n-k}^t(Ch_k^t(a,b),Ch_{k+1}^t(a,b))\\
Ch_n^t(a,b)&=aCh_n^t(1,0)+bCh_n^t(0,1)
\end{split}
\end{equation}
If we denote a particular Chebyshev polynomials $Ch_n^t(0,1)$ by $s_n(t)$, then from the recursion 
$$Ch_n^t(1,0)=Ch_{n-1}^t(0,-1)=-Ch_{n-1}^t(0,1)=-s_{n-1}(t),$$
 and arbitrary Chebyshev polynomials are expressed as  linear combinations of $s_{n-1}(t)$ and   $s_n(t)$ as follows :
 $$Ch_n^t(a,b)=-as_{n-1}(t)+bs_n(t)$$
The following two Chebyshev polynomials frequently appear along with $s_n(t)$:
\begin{equation*}
\begin{split}
v_n(t)&:=s_{n+1}(t)-s_n(t)=Ch_n^t(1,t-1)=Ch^t_{n+1}(1,1)\\
t_n(t)&:=s_{n+1}(t)-s_{n-1}(t)=Ch_n^t(2,t).\\
\end{split}
\end{equation*}
Note that $T_n(t)=\frac{1}{2}t_{n}(2t)$, $U_n(t)=s_{n+1}(2t)$, $V_n(t)=v_n(2t)$, and $W_n(t)=(-1)^nv_n(-2t)$ are exactly the classical Chebyshev polynomials of the first, second, third, and fourth  kind, respectively. 

\subsection{$\bf n$-Chebyshev polynomials}\label{n-Cheby}
\begin{definition}
For a given integer sequence ${\bf n}=(n_1,n_2,n_{3},\cdots)$,
we will call the  family of polynomials $\{g_n(t)\}$ defined by a three-term recursion
$$g_{k+1}(t)=n_{k} tg_k(t)-g_{k-1}(t)$$ 
${\bf n}$-$Chebyshev$ $polynomials$.
\end{definition}
We will denote the ${\bf n}$-Chebyshev polynomials $g_k(t)$ with the initial condition 
$g_0(t)=a, g_1(t)=b$ by $\mathcal Ch_k^{\bf n}(a,b)$. 
If we denote a particular ${\bf n}$-Chebyshev polynomials $\mathcal Ch_k^{\bf n}(0,1)$ by $s^{{\bf n}}_k(t)$, and
a sequence  $1$-shifted  to the right from ${\bf n}$ by $\sigma({\bf n})$, that is,  
$$\sigma^j({\bf n})=(n_{1+j},n_{2+j},\cdots)$$
then we have
$$\mathcal Ch_m^{\bf n}(a,b)=\mathcal Ch_{m-j}^{\sigma^j({\bf n})}(\mathcal Ch_j^{\bf n}(a,b),\mathcal Ch_{j+1}^{\bf n}(a,b))$$
and
$$\mathcal Ch_k^{\bf n}(a,b)=-as^{\sigma({\bf n})}_{k-1}(t)+bs^{{\bf n}}_k(t).$$
Note that $s^{{\bf n}}_{-1}(t) = -1$ by the recursive definition at $k=0$.
One may consider $S_k^\ns(t): = s_{k+1}^\ns(t)$, because $k$ coincides with the degree of $S^\ns_k(t)$, 
and note that only the first $k$ elements of $\ns$ are used in defining $s_{k+1}^\ns(t)$. 
We will continue to use $s_{k+1}^\ns(t)$ for consistency with our previous papers.

The following lemma can be easily proven using induction based on the definition of $s^{\bf n}_k$.
\begin{lemma}\label{n-cheby-1}
For an integer sequence ${\bf n}=(n_1,n_2,n_{3},\cdots)$, the following identities hold for $k=0,1,\cdots$.
\begin{enumerate}
\item[\rm (i)]  $s^{-\bf n}_{k}(t)=(-1)^{k+1}s^{\bf n}_{k}(t)$
\item [\rm (ii)]  $s^{-\bf n}_{k}(t)=s^{\bf n}_{k}(-t)$
\end{enumerate}
For ${\bf m} = (m_1, m_2, \cdots)$ and ${\bf l} = (l_1, l_2, \cdots)$ with $m_i = (-1)^i n_i$ and $l_i=(-1)^{i+1} n_i = -m_i$ $^\forall i \in \mathbb{N}$,
\begin{enumerate}
\item[\rm (iii)] $s^{\bf m}_{2k}(t) =-s^{\bf l}_{2k}(t)={\sqrt{-1}}s^{\bf n}_{2k}(\sqrt{-1}t)$
\item[\rm (iv)] $s^{\bf m}_{2k+1}(t) = s^{\bf l}_{2k+1}(t) = s^{\bf n}_{2k+1}(\sqrt{-1} t)$
\end{enumerate}

\end{lemma}
\begin{lemma}\label{ni-root-identity}
    For an integer sequence ${\bf n}=(n_1, n_2, \cdots)$ and $k \geq 1$,
    $$
    s^{\bf n}_{k+1}(t) = 
    \begin{cases}
        -s^{(n_3, n_4, \cdots, n_k)}_{k-1}(t) & \text{ if } n_1 = 0, \\
        -s^{(n_1,n_2,\cdots,n_{j-1}+n_{j+1},n_{j+2}, \cdots,n_k)}_{k-1}(t) & \text{ if } n_j = 0 \text{ for some } 1 < j < k, \\
        -s^{(n_1, n_2, \cdots, n_{k-2})}_{k-1}(t) & \text{ if } n_k = 0.
    \end{cases}
    $$
\end{lemma}

\begin{proof}
    For the cases $j=1$ or $j=k$, one can directly compute from the definition of ${\bf n}$-Chebyshev polynomials. If $k=j+1$, \begin{align*}
        s^{\bf n}_{j+2}(t) &= n_{j+1}t s^{\bf n}_{j+1}(t) - s^{\bf n}_{j}(t) \\
        &= n_{j+1}t (-s^{\bf n}_{j-1}(t)) - (n_{j-1}t s^{\bf n}_{j-1}(t) - s^{\bf n}_{j-2}(t)) \\
        &= -\left( (n_{j-1} + n_{j+1})ts^{\bf n}_{j-1}(t) - s^{\bf n}_{j-2}(t) \right) \\
        &= -s^{(n_1,n_2,\cdots,n_{j-1}+n_{j+1})}_{j} (t).
    \end{align*}
    Note that $s^{\bf n}_{j+1}(t) = -s^{\bf n}_{j-1}(t)$ when $n_{j}=0$. The conclusion for $k>j$ is deduced from the above calculation and the definition of $\bf{n}$-Chebyshev polynomials.   
\end{proof}

Note that $s^{\bf n}_{k}(t)$ is a generalized version of the classical Chebyshev polynomial of the second kind. And we introduce a generalized version of the third and fourth kind, $v^{\bf n}_{k}(t)$ and $v^{\bf n}_{k}(-t)$ respectively.

\begin{lemma}\label{n-cheby}
Let $v^{{\bf n}}_k(t)$ be the ${\bf n}$-Chebyshev polynomials defined by
$$
v^{{\bf n}}_k(t):=\mathcal Ch_{k+1}^{\bf n}(1,1).\\
$$
Then
\begin{enumerate}
\item[\rm (i)]  $v^{\bf n}_{k+2}(t)=n_{k+2}t v^{\bf n}_{k+1}(t)-v^{\bf n}_{k}(t)$, that is, $v^{\bf n}_k(t)=\mathcal Ch_{k}^{\sigma({\bf n})}(1,n_1t-1)$
\item [\rm (ii)]
$v^{\bf n}_k(t)=s^{\bf n}_{k+1}(t)-s^{\sigma({\bf n})}_{k}(t)$
\item [\rm (iii)]  $v^{\bf n}_k(-t)=(-1)^k(s^{\bf n}_{k+1}(t)+s^{\sigma({\bf n})}_{k}(t))=(-1)^n\mathcal Ch_{k+1}^{\bf n}(-1,1)$
\item [\rm (iv)] $v^{-\bf n}_{k}(t)=v^{\bf n}_{k}(-t)$
\end{enumerate}
\end{lemma}
\begin{proof}
    The statements can be easily proved by induction, based on the definitions of $v^{\bf n}_k$ and $s^{\bf n}_k$. We omit the proofs.
\end{proof}
Now, let's consider the following three polynomials associated with the classical Chebyshev polynomial of the first kind:

\begin{itemize}
\item $T^{\bf n}_{k}(t):=s^{\bf n}_{k+1}(t) -s^{\bf n}_{k-1}(t)=n_{k} t s^{\bf n}_{k}(t) -2s^{\bf n}_{k-1}(t)$
\item $\tilde{T}^{\bf n}_{k}(t):=s^{\bf n}_{k+1}(t) -s^{\sigma({\bf n})}_{k-1}(t)=n_{k} t s^{\bf n}_{k}(t) -s^{\bf n}_{k-1}(t)-s^{\sigma({\bf n})}_{k-1}(t)$
\item $\hat{T}^{\bf n}_{k}(t):=s^{\bf n}_{k+1}(t)+s^{\sigma({\bf n})}_{k-1}(t)=n_{k} t s^{\bf n}_{k}(t) -s^{\bf n}_{k-1}(t)+s^{\sigma({\bf n})}_{k-1}(t)$
\end{itemize}
From the definition, we have the following identity immediately.
\begin{lemma}\label{T-s identity}
For any integer squence $(n_1,n_2,\cdots,n_{k})$ we have  
$$s^{(n_1,n_2,\cdots,n_{k})}_{k+1}(t)=\frac{1}{2}T_{k}^{(n_1,n_2,\cdots,n_{k-1}, 2n_{k})}(t).$$
\end{lemma}
\begin{proof}
The identity is proved as follows:
\begin{equation*}
\begin{split}
T_{k}^{(n_1,n_2,\cdots,n_{k-1}, 2n_{k})}(t)&=s_{k+1}^{(n_1,n_2,\cdots,n_{k-1}, 2n_{k})}(t)-s_{k-1}^{(n_1,n_2,\cdots,n_{k-1}, 2n_{k})}(t)\\
&=2n_{k}ts_{k}^{(n_1,n_2,\cdots,n_{k-1})}(t)-s_{k-1}^{(n_1,n_2,\cdots,n_{k-1})}(t)-s_{k-1}^{(n_1,n_2,\cdots,n_{k-1}, 2n_{k})}(t)\\
&=2(n_{k}ts_{k}^{(n_1,n_2,\cdots,n_{k-1})}(t)-s_{k-1}^{(n_1,n_2,\cdots,n_{k-1})}(t))\\
&=2s^{(n_1,n_2,\cdots,n_{k})}_{k+1}(t).
\end{split}
\end{equation*}
\end{proof}

\subsection{Continuant}
From the three term recursive definition of n-Chebyshev polynomial, we are led to consider a continuant, the determinant of a tridiagonal matrix, which in turn reminds us a connection with a continued fraction. We recall some basic properties of the continuant, and then study some divisibility properties for later use in this subsection.
\begin{definition}
Each of the following two kinds of tridiagonal determinant is called the $continuant$.
$$
K_n(c_1,\cdots,c_n):=\det 
\begin{pmatrix}
c_1 & 1 & & &\\
1&c_2 &1 & &\\
& \ddots & \ddots & \ddots &\\
&& 1&c_{n-1}&1\\
&&&1&c_n
 \end{pmatrix}
$$
$$
K^+_n(c_1,\cdots,c_n):=\det 
\begin{pmatrix}
c_1 & 1 & & &\\
-1&c_2 &1 & &\\
& \ddots & \ddots & \ddots &\\
&& -1&c_{n-1}&1\\
&&&-1&c_n
 \end{pmatrix}
$$
\end{definition}
We set $K_0:=1, K_{-1}:=0$ and  $K_0^+:=1, K_{-1}^+:=0$. Then we have the recursion for $K$ and $K^+$:
$$
K_n(c_1,\cdots,c_n)=c_nK_{n-1}(c_1,\cdots,c_{n-1})-K_{n-2}(c_1,\cdots,c_{n-2})$$
and
 $$
K^+_n(c_1,\cdots,c_n)=c_nK^+_{n-1}(c_1,\cdots,c_{n-1})+K^+_{n-2}(c_1,\cdots,c_{n-2}).$$ 
Using these recursion  we get easily the following lemma by induction, and we omit the proof.
\begin{lemma}\label{continuant-lemma1} We have
\begin{enumerate}[label=\normalfont(\roman*)]
\item
$
K_{2n}(c_1x,c_2,c_3x\cdots,c_{2n-1}x,c_{2n})=K_{2n}(c_1,c_2x,c_3\cdots,c_{2n-1},c_{2n}x)
$
\\
$
K^+_{2n}(c_1x,c_2,c_3x\cdots,c_{2n-1}x,c_{2n})=K^+_{2n}(c_1,c_2x,c_3\cdots,c_{2n-1},c_{2n}x)
$ 
\item
$
K_{2n+1}(c_1x,c_2,c_3x\cdots,c_{2n-1}x,c_{2n},c_{2n+1}x )=xK_{2n+1}(c_1,c_2x,c_3\cdots,c_{2n-1},c_{2n}x, c_{2n+1})$
\\
$K^+_{2n+1}(c_1x,c_2,c_3x\cdots,c_{2n-1}x,c_{2n},c_{2n+1}x )=xK^+_{2n+1}(c_1,c_2x,c_3\cdots,c_{2n-1},c_{2n}x, c_{2n+1})$
\item \label{continuant-lemma1 K reversible}
$K_n(c_1,\cdots,c_n)=K_n(c_n,\cdots,c_{1})$
\\
$K^+_n(c_1,\cdots,c_n)=K^+_n(c_n,\cdots,c_{1})$
\item
$K_n(-c_1,\cdots,-c_n)=(-1)^nK_n(c_1,\cdots,c_n)$
\\
$K^+_n(-c_1,\cdots,-c_n)=(-1)^nK^+_n(c_1,\cdots,c_{n})$
\end{enumerate}
\end{lemma}
The following proposition shows some of the identities between the two continuants, which will be used later.
\begin{proposition}\label{Kplus&K}
Let $x=-t^2\in \bc$. Then we have
\begin{enumerate}
\item[\rm (i)]
$
K^+_{2n}(c_1,c_2x,c_3,c_4x\cdots,c_{2n-1},c_{2n}x)=(-1)^nK_{2n}(c_1t,c_2t,\cdots,c_{2n}t)$ 
\\ 
$
K^+_{2n-1}(c_1, c_2x,c_3,c_4x\cdots,c_{2n-1})=\frac{(-1)^{n-1}}{t}K_{2n-1}(c_1t,c_2t,c_3t\cdots,c_{2n-1}t)$ 
\item[\rm (ii)]
$
K^+_{2n-1}(c_1x,c_2\cdots,c_{2n-2},c_{2n-1}x )=(-1)^ntK_{2n-1}(c_1t,c_2t\cdots,c_{2n-1}t)$
\\ 
$K^+_{2n}(c_1x,c_2\cdots,c_{2n-1}x,c_{2n} )=(-1)^{n}K_{2n}(c_1t,c_2t\cdots,c_{2n}t)$
\end{enumerate}
\end{proposition}
\begin{proof}
We only provide a proof for (i) here, since the similar argument can be applied to prove (ii).
For $n=1$, the identities follow from
$$K^+_{1}(c_1)=c_1,\,\,K^+_{2}(c_1,c_2x)=\begin{vmatrix}
c_1 & 1\\
-1 & c_2x 
 \end{vmatrix}=c_1c_2x+1=-c_1c_2t^2+1$$
 and
 $$K_{1}(c_1t)=c_1t, \,\,K_{2}(c_1t,c_2t)=\begin{vmatrix}
c_1t & 1\\
1 & c_2t 
 \end{vmatrix}=c_1c_2t^2-1.$$
 To proceed with induction, we assume that the identities hold for all $n\leq k$.  Then we have
 \begin{equation*}
\begin{split}
K^+_{2k+1}&(c_1,c_2x,c_3,c_4x\cdots,c_{2k+1})\\
&=c_{2k+1}K^+_{2k}(c_1,c_2x,c_3,c_4x\cdots,c_{2k}x)+K^+_{2k-1}(c_1,c_2x,c_3,c_4x\cdots,c_{2k-1})\\
&=(-1)^kc_{2k+1}K_{2k}(c_1t,c_2t,\cdots,c_{2k}t)+\frac{(-1)^{k-1}}{t}K_{2k-1}(c_1t,c_2t,c_3t\cdots,c_{2k-1}t)\\
&=(-1)^k\frac{1}{t}(c_{2k+1}tK_{2k}(c_1t,c_2t,\cdots,c_{2k}t)-K_{2k-1}(c_1t,c_2t,c_3t\cdots,c_{2k-1}t))\\
&=\frac{(-1)^{k}}{t}K_{2k+1}(c_1t,c_2t,c_3t\cdots,c_{2k+1}t)
\end{split}
\end{equation*}
and
 \begin{equation*}
\begin{split}
K^+_{2k+2}&(c_1,c_2x,c_3,c_4x\cdots,c_{2k+1},c_{2k+2}x)\\
&=c_{2k+2}xK^+_{2k+1}(c_1,c_2x,c_3,c_4x\cdots,c_{2k+1})+K^+_{2k}(c_1,c_2x,c_3,c_4x\cdots,c_{2k-1},c_{2k}x)\\
&=(-c_{2k+2}t^2)\frac{(-1)^{k}}{t}K_{2k+1}(c_1t,c_2t,c_3t\cdots,c_{2k+1}t)+(-1)^kK_{2k}(c_1t,c_2t,\cdots,c_{2k}t)\\
&=(-1)^{k+1}\left (c_{2k+2}tK_{2k+1}(c_1t,c_2t,c_3t\cdots,c_{2k+1}t)-K_{2k}(c_1t,c_2t,\cdots,c_{2k}t)\right)\\
&=(-1)^{k+1}K_{2k+1}(c_1t,c_2t,c_3t\cdots,c_{2k+2}t),
\end{split}
\end{equation*}
which implies that (i) is also true for $n=k+1$ as desired. 
This completes the proof.
\end{proof}

The following proposition is well-known.
\begin{proposition}[\cite{CO}, \cite{MO}]\label{KN}
Let 
$$M(c_1,\cdots,c_n):=
\begin{pmatrix}
c_1& -1\\
1& 0 
 \end{pmatrix}
\begin{pmatrix}
c_2& -1\\
1& 0 
 \end{pmatrix}
\cdots
\begin{pmatrix}
c_n& -1\\
1& 0 
 \end{pmatrix}
$$
and
$$M^+(c_1,\cdots,c_n):=
\begin{pmatrix}
c_1& 1\\
1& 0 
 \end{pmatrix}
\begin{pmatrix}
c_2& 1\\
1& 0 
 \end{pmatrix}
\cdots
\begin{pmatrix}
c_n& 1\\
1& 0 
 \end{pmatrix}.
$$
Then 
$$M(c_1,\cdots,c_n)=\begin{pmatrix}
K_n(c_1,\cdots,c_n)& -K_{n-1}(c_1,\cdots,c_{n-1})\\
K_{n-1}(c_2,\cdots,c_{n})& -K_{n-2}(c_2,\cdots,c_{n-1}) 
 \end{pmatrix}
$$
and
$$M^+(c_1,\cdots,c_n)=\begin{pmatrix}
K^+_n(c_1,\cdots,c_n)& K^+_{n-1}(c_1,\cdots,c_{n-1})\\
K^+_{n-1}(c_2,\cdots,c_{n})& K^+_{n-2}(c_2,\cdots,c_{n-1}) 
 \end{pmatrix}.
$$
\end{proposition}
The following is known as Euler's Identity and play a crucial role for the combinatorial properties of the continuants. (See \cite{GKP}, \cite{MO}, or \cite{ustinov2006short} for references.)
\begin{proposition}[Euler Identity]\label{Euler Identity}
$K_{k-i-1}(c_{i+1},\cdots,c_{k-1})K_{l-j-1}(c_{j+1},\cdots,c_{l-1})$ is equal to 
$$K_{j-i-1}(c_{i+1},\cdots,c_{j-1})K_{l-k-1}(c_{k+1},\cdots,c_{l-1})+K_{k-j-1}(c_{j+1},\cdots,c_{k-1})K_{l-i-1}(c_{i+1},\cdots,c_{l-1}).$$
\end{proposition} 
Applying  Proposition \ref{Euler Identity} to the case of $k=j+1$, we get
\begin{corollary}\label{EI-cor}
$K_{j-i}(c_{i+1},\cdots,c_{j})K_{l-j-1}(c_{j+1},\cdots,c_{l-1})$ is equal to 
$$K_{j-i-1}(c_{i+1},\cdots,c_{j-1})K_{l-j-2}(c_{j+2},\cdots,c_{l-1})+K_{l-i-1}(c_{i+1},\cdots,c_{l-1}).$$
\end{corollary}
Using Corollary \ref{EI-cor} we get the following proposition.
\begin{proposition}[\cite{Jo-Kim2}]\label{EI-1}
For any integer sequence ${\bf c}=(c_1,\cdots,c_n)$, ${\overleftarrow{\bf c}} = (c_n, \cdots, c_1)$, and nonzero $x$, the following hold.
\begin{enumerate}
\item[\rm (i)] 
$K_{2n+1}({\bf c}, x,{\bf c})=K_n({\bf c})\left(K_{n+1}({\bf c},x)-K_{n-1}(c_2,\cdots,c_{n})\right )$
\\
$K_{2n+1}({\bf c}, x,-{\bf c})=(-1)^nK_n({\bf c})\left(K_{n+1}({\bf c},x)+K_{n-1}(c_2,\cdots,c_{n})\right )$
\\
$K_{2n+1}({\bf c}, x,\overleftarrow{\bf c})=K_n({\bf c})\left (K_{n+1}({\bf c},x)-K_{n-1}(c_1,\cdots,c_{n-1})\right )$
\\
$K_{2n+1}({\bf c}, x,-\overleftarrow{\bf c})=(-1)^nxK_n({\bf c})^2$
\item[\rm (ii)]
$K_{2n}(\overleftarrow{\bf c},{\bf c})=K_{n}({\bf c})^2-K_{n-1}(c_2,\cdots,c_n)^2$
\\
 $K_{2n}({\bf c},{\bf c})
=K_{n}( {\bf c})\left( K_{n}( {\bf c})-K_{n-2}( c_2,\cdots,c_{n-1})\right)-1
$
\\
$K_{2n}({\bf c},-{\bf c})
=(-1)^n\left(K_{n}( {\bf c})\left( K_{n}( {\bf c})+K_{n-2}( c_2,\cdots,c_{n-1})\right)+1\right)
$
\end{enumerate}
\end{proposition}
For each ${\bf c}=(c_1,\cdots,c_n)$, consider the following three polynomials in $x$.
\begin{itemize}
\item $q_{\bf c}(x):=xK_n(c_1,\cdots,c_n)-2K_{n-1}(c_1,\cdots,c_{n-1})$
\item $p_{\bf c}(x):=xK_n(c_1,\cdots,c_n)-K_{n-1}(c_1,\cdots,c_{n-1})-K_{n-1}(c_2,\cdots,c_n)$
\item $w_{\bf c}(x):=xK_n(c_1,\cdots,c_n)-K_{n-1}(c_1,\cdots,c_{n-1})+K_{n-1}(c_2,\cdots,c_n)$
\end{itemize}
Then the first three identities in (i) of Proposition \ref{EI-1} can be reformulated as follows.
\begin{proposition}\label{p-q}
Let ${\bf c}=(c_1,\cdots,c_n)$.
Then
\begin{equation*}
    \begin{split}    K_{2n+1}(c_1,\cdots,c_n,x,c_n,\cdots,c_1)&=K_n(c_1,\cdots,c_n)q_{\bf c}(x)\\     K_{2n+1}(c_1,\cdots,c_n,x,c_1,\cdots,c_n)&=K_n(c_1,\cdots,c_n)p_{\bf c}(x)\\
    K_{2n+1}(c_1,\cdots,c_n,x,-c_1,\cdots,-c_n)&=(-1)^nK_n(c_1,\cdots,c_n)w_{\bf c}(x).
    \end{split}
\end{equation*}
\end{proposition}
More generally, we get the following.
\begin{proposition}\label{cont-fact}
Let ${\bf c}=(c_1,\cdots,c_n)$. 
Then
\begin{enumerate}
\item[\rm (i)]
$K_{kn+k+n}({\bf c},x_1,{\bf c},x_2,\cdots,{\bf c},x_k, {\bf c})=K_n({\bf c})K_k(p_{\bf c}(x_1),\cdots,p_{\bf c}(x_k))$
\item[\rm (ii)]   
$K_{2kn+2k+n}({\bf c},x_1,{\overleftarrow{\bf c}},x_2,{\bf c},x_3,{\overleftarrow{\bf c}},\cdots, x_{2k}, {\bf c})=K_n({\bf c})K_{2k}(q_{\bf c}(x_1),q_{{\overleftarrow{\bf c}}}(x_2),\cdots,q_{\overleftarrow{\bf c}}(x_{2k}))$   
\\
$K_{2kn+2k+2n+1}({\bf c},x_1,{\overleftarrow{\bf c}},x_2,{\bf c},\cdots, x_{2k+1},{\overleftarrow{\bf c}})=K_n({\bf c})K_{2k+1}(q_{\bf c}(x_1),q_{\overleftarrow{\bf c}}(x_2),\cdots,q_{\bf c}(x_{2k+1}))$   
\item[\rm (iii)]   
 $K_{2kn+2k+n}({\bf c},x_1,-{\bf c},x_2,{\bf c},\cdots,-{\bf c},x_{2k}, {\bf c})=(-1)^{kn}K_n({\bf c})K_k(w_{\bf c}(x_1),w_{\overleftarrow{\bf c}}(x_2),\cdots,w_{\overleftarrow{\bf c}}(x_{2k}))$
\\
 $K_{2kn+2k-1}({\bf c},x_1,-{\bf c},x_2,{\bf c},\cdots,x_{2k-1}, -{\bf c})=(-1)^{kn}K_n({\bf c})K_k(w_{\bf c}(x_1),w_{\overleftarrow{\bf c}}(x_2),\cdots,w_{\bf c}(x_{2k-1}))$
\end{enumerate}
\end{proposition}
\begin{proof}
Since the proofs for (i) through (iii) follow a similar pattern, we will provide a detailed proof for (i).
For $k=1$, the identity 
$$K_{2n+1}({\bf c},x,{\bf c})=K_n({\bf c})K_1(p_{\bf c}(x))$$
follows from (i) of Proposition \ref{EI-1}.
We proceed  by induction on $k$.  Assume that the identity holds for all $k < i$. Then 
$$K_{jn+j+n}({\bf c},x_1,{\bf c},x_2,\cdots,x_{j},{\bf c})=K_n({\bf c})K_{j}(p_{\bf c}(x_1),\cdots,p_{\bf c}(x_{j}))$$
for all $j < i$. Since
\begin{equation*}
\begin{split}
K_{ni+i-1}&({\bf c},x_1,{\bf c},x_2,\cdots,x_{i-1},{\bf c})K_{2n+1}({\bf c},x_{i},{\bf c})\\
&=K_n({\bf c})K_{ni+i+n}({\bf c},x_1,{\bf c},\cdots,x_i,{\bf c})+K_n({\bf c})K_{ni+i-n-2}({\bf c},x_1,{\bf c},\cdots,x_{i-2},{\bf c})
\end{split}
\end{equation*}
from Euler identity,  we have, by induction hypothesis, 
\begin{equation*}
\begin{split}
K_n({\bf c})^2&K_{i-1}(p_{\bf c}(x_1),\cdots,p_{\bf c}(x_{i-1}))K_1(p_{\bf c}(x_i))\\
&=K_n({\bf c})K_{ni+i+n}({\bf c},x_1,{\bf c},\cdots,x_i,{\bf c})+K_n({\bf c})K_{i-2}(p_{\bf c}(x_1),\cdots,p_{\bf c}(x_{i-2}))
\end{split}
\end{equation*}
and thus
\begin{equation*}
\begin{split}
K&_{ni+i+n}({\bf c},x_1,{\bf c},\cdots,x_i,{\bf c})\\
&=K_n\left({\bf c})(K_{i-1}(p_{\bf c}(x_1),\cdots,p_{\bf c}(x_{i-1}))K_1(p_{\bf c}(x_i))-K_{i-2}(p_{\bf c}(x_1),\cdots,p_{\bf c}(x_{i-2}))\right)\\
&=K_n\left({\bf c})p_{\bf c}(x_i)K_{i-1}(p_{\bf c}(x_1),\cdots,p_{\bf c}(x_{i-1}))-K_{i-2}(p_{\bf c}(x_1),\cdots,p_{\bf c}(x_{i-2}))\right)\\
&=K_n({\bf c})K_{i}(p_{\bf c}(x_1),\cdots,p_{\bf c}(x_{i})).
\end{split}
\end{equation*}
which implies that the identity holds for $k=i$ also.
\end{proof}

\subsection{Relationship between Continuants and $\bf n$-Chebyshev polynomials}
The continuant is related to ${\bf n}$-Chebyshev polynomials by letting $c_i=n_it$. The following lemma  illustrates this relationship.
\begin{lemma}\label{s-v-continuant}
Let ${\bf n}=(n_1,n_2,n_{3},\cdots)$. Then
\begin{enumerate}[label=\normalfont(\roman*)]
\item \label{s-v-continuant s equal to K}
$s^{\bf n}_{m+1}(t)=K_{m}( n_1t, n_2t,\cdots,n_mt)$
\item
$v^{\bf n}_m(t)=K_{m+1}(1, n_1t, n_2t,\cdots,n_mt)=K_{m}( n_1t, n_2t,\cdots,n_mt)-K_{m-1}(n_2t,\cdots,n_mt)$
\item
$v^{\bf n}_m(-t)=(-1)^n\left (K_{m}( n_1t,\cdots,n_{m}t)+K_{m-1}(n_2t,\cdots,n_mt)\right )$
\begin{equation*}
\begin{split}
{\rm  (iv)}\,\,  n_mt&K_{m+1}(1, n_1t,\cdots,n_{m-1}t)\\
&=K_{m}( n_1t,\cdots,n_mt)-K_{m-1}(n_2t,\cdots,n_mt)+K_{m-1}(1, n_1t,\cdots,n_{m-2}t)\qquad
\end{split}
\end{equation*}
\end{enumerate}
\end{lemma}
\begin{proof}
For $m=1$, the identity (i) and the first equality of (ii) follow from 
$$s^{\bf n}_1(t)=1=K_0, \quad s^{\bf n}_2(t)=n_1t=K_1(n_1t)$$  
and
$$v^{\bf n}_0(t)=1=K_1(1), \quad v^{\bf n}_1(t)=n_1t-1=K_2(1, n_1t)$$ 
 To proceed with induction, we assume that the identities hold for all $m \leq k$.  Then we have
\begin{equation*}
\begin{split}
s^{\bf n}_{k+1}(t)&= n_{k}ts^{\bf n}_{k}(t)-s^{\bf n}_{k-1}(t)\\
&=n_{k}tK_{k-1}(n_1t,\cdots,n_{k-1}t)-K_{k-1}(n_1t,\cdots,n_{k-2}t)\\
&=K_{k+1}(n_1t,\cdots,n_{k+1}t)
\end{split}
\end{equation*}
and 
\begin{equation*}
\begin{split}
v^{\bf n}_{k+1}(t)&= n_{k+1}tv^{\bf n}_{k}(t)-v^{\bf n}_{k-1}(t)\\
&= n_{k+1}t  K_{k+1}(1, n_1t, n_2t,\cdots,n_kt)- K_{k}(1, n_1t, n_2t,\cdots,n_{k-1}t)\\
&=K_{k+2}(1, n_1t, n_2t,\cdots,n_{k+1}t)
\end{split}
\end{equation*}
which implies that (i) and the first equality of (ii) are also true for $m=k+1$ as desired. 

The second equality of (ii) follows immediately from Lemma \ref{n-cheby}, (ii).
(iii) is proved as follows.
\begin{equation*}
\begin{split}
v^{\bf n}_m(-t)&=K_{m}( -n_1t, -n_2t,\cdots,-n_mt)-K_{m-1}(-n_2t,\cdots,-n_mt) \\
&=(-1)^m K_{m}( n_1t, n_2t,\cdots,n_mt)-(-1)^{m-1}K_{m-1}(n_2t,\cdots,n_mt) \\
&=(-1)^m \left (K_{m}( n_1t,\cdots,n_{m}t)+K_{m-1}(n_2t,\cdots,n_mt)\right )
\end{split}
\end{equation*}

(iv) follows immediately from the fact that 
$$K_{m+1}(1, n_1t, n_2t,\cdots,n_mt)=n_mtK_{m}(1, n_1t,\cdots,n_{m-1}t)-K_{m-1}(1, n_1t,\cdots,n_{m-2}t).$$
\end{proof}


Applying Proposition \ref{EI-1} to the case where $(c_1,\cdots,c_k)$ is $(n_1t,\cdots,n_kt)$ and $x$ is $n_{k+1}$, along with the relationship between the continuants and ${\bf n}$-Chebyshev polynomials from 
Lemma \ref{s-v-continuant}, we obtain the following identities.
\begin{proposition}\label{factor-cont-prop}
Let  ${\bf n}=(n_1,n_2,\cdots, n_{k})$ be an integer sequence. Then the following holds.
\begin{enumerate}
\item[\rm (i)] $s^{({\bf n}, n_{k+1}, {\overleftarrow{\bf n}})}_{2k+2}(t)=s^{\bf n}_{k+1}(t)T^{({\bf n},n_{k+1})}_{k+1}(t)$
\item[\rm (ii)] $s^{({\bf n},n_{k+1},{\bf n})}_{2k+2}(t)=s^{\bf n}_{k+1}(t)\tilde{T}^{({\bf n},n_{k+1})}_{k+1}(t)$
\item[\rm (iii)] $s^{({\bf n},n_{k+1},-{\bf n})}_{2k+2}(t)=s^{\bf n}_{k+1}(t)\hat{T}^{({\bf n},n_{k+1})}_{k+1}(t)$
\item[\rm(iv)]$s^{({\overleftarrow{\bf n}},{\bf n})}_{2k+1}(t)=(-1)^kv^{\bf n}_k(t)v^{\bf n}_k(-t)$
\item[\rm (v)] $s^{({\bf n},{\bf n})}_{2k+1}(t)=s^{\bf n}_{k+1}(t)\tilde{T}^{\bf n}_{k}(t)-1$
\item[\rm (vi)] $s^{({\bf n},-{\bf n})}_{2k+1}(t)=(-1)^k\big(s^{\bf n}_{k+1}(t)\hat{T}^{\bf n}_{k}(t)+1\big)$
\end{enumerate}
\end{proposition}
\begin{proof}
The identities (i), (ii), and (iii) follow from (i) of Proposition \ref{EI-1}: For example, 
\begin{equation*}
\begin{split}
s^{({\bf n}, n_{k+1}, {\overleftarrow{\bf n}})}_{2k+2}(t)
&=K_{2k+1}(n_1t,\cdots,n_kt,n_{k+1}t,n_kt,\cdots,n_1t)\\
&=K_k(n_1t,\cdots,n_kt)\left(K_k(n_1t,\cdots,n_{k+1}t)-K_{k-1}(n_1t,\cdots,n_{k-1}t)\right)\\
&=s^{\bf n}_{k+1}(t)\left(s^{\bf n}_{k+2}(t) -s^{\bf n}_{k}(t)\right)\\
&=s^{\bf n}_{k+1}(t)T^{\bf n}_{k+1}(t)
\end{split}
\end{equation*}
and similarly for (ii) and (iii). 
The identity (iv) follows from (ii) of Proposition \ref{EI-1} and Lemma \ref{n-cheby}:
\begin{equation*}
\begin{split}
s^{({\overleftarrow{\bf n}},{\bf n})}_{2k+1}(t)&=K_{2k}( n_kt,\cdots,n_1t, n_1t,\cdots,n_{k}t)\\
&=K_{k}(n_1t,\cdots,n_{k}t)^2-K_{k-1}(n_2t,\cdots,n_{k}t)^2\\
&=s^{{\bf n}}_{k+1}(t)^2-s^{\sigma({\bf n})}_{k}(t)^2\\
&=(s^{{\bf n}}_{k+1}(t)+s^{\sigma({\bf n})}_{k}(t))(s^{{\bf n}}_{k+1}(t)-s^{\sigma({\bf n})}_{k}(t)))\\
&=(-1)^kv^{\bf n}_k(t)v^{\bf n}_k(-t)
\end{split}
\end{equation*}
The identities (v) and (vi) are also proved similarly, using (ii) of Proposition \ref{EI-1}.
\end{proof}

Further, the relationships given in Lemma \ref{s-v-continuant} can be used to derive more identities for ${\bf n}$-Chebyshev polynomials. Let us introduce notations that help simplify statements and arguments.
For an integer sequence ${\bf n}=(n_1, n_2, \cdots, n_k)$ and fixed $1 \leq i \leq k$, let
\begin{align*}
{\bf n}_L := (n_1, n_2, \cdots, n_{i-1}), &\, {\bf n}_{L'} := (n_1, n_2, \cdots, n_{i-2}) \\
{\bf n}_R := (n_{i+1}, \cdots, n_{k-1}, n_{k}), &\,{\bf n}_{R'} := (n_{i+2}, \cdots, n_{k-1}, n_{k})
\end{align*}
be subsequences of ${\bf n}$,
and
\begin{align*}
    s^{{\bf n}_L}(t) :=s^{(n_1, n_2, \cdots, n_{i-1})}_i(t), &\, s^{{\bf n}_{L'}}(t) := s^{(n_1, n_2, \cdots, n_{i-2})}_{i-1}(t) \\
    s^{{\bf n}_R}(t) := s^{(n_{i+1}, n_{i+2}, \cdots, n_{k})}_{k-i+1}(t), &\, s^{{\bf n}_{R'}}(t) := s^{(n_{i+2}, n_{i+3}, \cdots, n_{k})}_{k-i}(t).
\end{align*} Note that above polynomials are well-defined when $i=1$ and $i=k$ since $s^{{\bf n}_*}_0(t) = 0$ and $s^{{\bf n}_*}_1(t) = 1$ are independent for all ${\bf n}_* \in \{{\bf n}_L,{\bf n}_{L'},{\bf n}_R,{\bf n}_{R'}\}$. 
When the length of ${\bf n}$ is $k$, all the integers in ${\bf n}$ are used to determine $s^{\bf n}_{k+1}(t)$, and hence we may simply denote $s^{\bf n}_{k+1}(t)$ by $s^{\bf n}$, and also $T^{{\bf n}}, \tilde{T}^{{\bf n}}$ and $\hat{T}^{{\bf n}}$ will be used instead of $T^{{\bf n}}_k(t), \tilde{T}^{{\bf n}}_k(t)$ and $\hat{T}^{{\bf n}}_k(t)$, respectively.



${\bf n}$-Chebyshev polynomials are all first order polynomials in each $n_i$, when they are considered as a function of $n_i$ with $t$ fixed.  The following proposition shows us that their coefficients are expressed in terms of $s^{\bf m}(t)$'s in such considerations.

\begin{proposition}\label{split s}
Let ${\bf n} = (n_1, n_2, \cdots, n_k)$ be an integer sequence. With the notations above, we have the following, for any $1 \leq i \leq k$:
\begin{enumerate}[label=\normalfont(\roman*)]
\item \label{split s 2terms}
$s^{\bf n}_{k+1}(t) = 
    n_i t s^{{\bf n}_L}_{i-1}(t) s^{{\bf n}_R}_{k+1-i}(t) - s^{{\bf n}_{L'}}_{i-2}(t) s^{{\bf n}_R}_{k+1-i}(t) - s^{{\bf n}_L}_{i-1}(t) s^{{\bf n}_{R'}}_{k-i}(t)
$
\item \label{split s 1term}
$s^{\bf n}_{k+1}(t) =
    n_i t s^{{\bf n}_L}_{i-1}(t) s^{{\bf n}_R}_{k+1-i}(t)-s^{({\bf n}_{L'}, n_{i-1} + n_{i+1}, {\bf n}_{R'})}_{k-1}(t)
$
\item \label{split s split ttilde}
$\tilde{T}^{{\bf n}}_{k}(t)
    = n_i t s^{
        ({\bf n}_R,{\bf n}_L)
    }_k(t)
    - s_{k-1}^{({\bf n}_{R}, {\bf n}_{L'})}(t) - s^{({{\bf n}_{R'}}, {\bf n}_{L})}_{k-1}(t)
$
\item \label{split s split that}
$\hat{T}^{{\bf n}}_{k}(t)
    = (-1)^{i-1}(n_i t s^{
    ({\bf n}_{R}, -{\bf n}_L)
    }_k(t)
    + s^{({\bf n}_{R}, -{\bf n}_{L'})}_{k-1}(t) - s^{({\bf n}_{R'}, -{\bf n}_L)}_{k-1}(t)
    )
$
\end{enumerate}
\end{proposition}

\begin{proof}
Since the statements are obvious when $i$ is either $1$ or $k$ by the definition of each polynomial (recall that $s^{{\bf n}}_{-1}(t)=-1$), we may assume that $1<i<k$. For the first two statements, 
by substituting $i=0, j=i, k=i+1, l=k+1$ and $c_m = n_m t$ for all $m$ in the Euler Identity (Proposition \ref{Euler Identity}), 
$K_{k}(n_1t,\cdots, n_kt)$ is equal to 
$$
K_{i}(n_1t,\cdots,n_it)K_{k-i}(n_{i+1}t,\cdots,n_kt)-K_{i-1}(n_1t,\cdots,n_{i-1}t)K_{k-i-1}(n_{i+2}t,\cdots,n_kt).
$$
Hence by Lemma \ref{s-v-continuant} we have
    \begin{align*}
    s^{\bf n} &= s^{({\bf n}_L, n_i)} s^{{\bf n}_R} - s^{{\bf n}_L} s^{{\bf n}_{R'}} \\
    &= n_i t s^{{\bf n}_L} s^{{\bf n}_{R}}-s^{{\bf n}_{L'}} s^{{\bf n}_R} - s^{{\bf n}_{L}} s^{{\bf n}_{R'}},
    \end{align*}
which proves that \ref{split s 2terms}, and by direct implication of Lemma \ref{ni-root-identity}, \ref{split s 1term} also holds true.

For (iii), recall that $\tilde{T}^{{\bf n}}_k = s^{{\bf n}}_{k+1} - s_{k-1}^{(n_2, \cdots, n_{k-1})}$. Notice that by (i),
  \begin{align*}
    \tilde{T}_k^{{\bf n}} &= n_i t \cdot P(i-1, i+1) - P(i-2, i+1) - P(i-1, i+2)
  \end{align*}
  where $P(j, l) := s^{(n_1, \cdots, n_j)}_{j+1} s^{(n_l, \cdots, n_k)}_{k-l+2} - s^{(n_2, \cdots, n_j)}_{j} s^{(n_l, \cdots, n_{k-1})}_{k-l+1}$.

  Observe that
  \begin{align*}
    P(j, l) & = K(n_1 t, \cdots, n_j t) K(n_l t, \cdots, n_k t) - K(n_2 t, \cdots, n_j t) K(n_l t, \cdots, n_{k-1} t) \\
    & = K(n_j t, \cdots, n_1 t) K(n_k t, \cdots, n_l t) - K(n_j t, \cdots, n_2 t) K(n_{k-1} t, \cdots, n_l t) \\
    & = K(n_j t, \cdots, n_1 t, n_k t, \cdots, n_l t)
  \end{align*}
  follows from Corollary \ref{EI-cor}. Therefore,
  \begin{align*}
    \tilde{T}^{{\bf n}}_k
    & = n_i t s^{(n_{i-1}, \cdots, n_1, n_k, \cdots, n_{i+1})}_{k}
    - s^{(n_{i-2}, \cdots, n_1, n_k, \cdots, n_{i+1})}_{k-1}
    - s^{(n_{i-1}, \cdots, n_1, n_k, \cdots, n_{i+2})}_{k-1} \\
    & = n_i t s^{(n_{i+1}, \cdots, n_k, n_1, \cdots, n_{i-1})}_{k}
    - s^{(n_{i+1}, \cdots, n_k, n_1, \cdots, n_{i-2})}_{k-1}
    - s^{(n_{i+2}, \cdots, n_k, n_1, \cdots, n_{i-1})}_{k-1},
  \end{align*}
  by the symmetry of the continuant (v) in Lemma \ref{continuant-lemma1} as desired.

Similarly for $\hat{T}^{{\bf n}}_k = s^{{\bf n}}_{k+1} + s_{k-1}^{(n_2, \cdots, n_{k-1})}$, one has
  \[
    \hat{T}_k^{\bf n} = n_i t \cdot Q(i-1, i+1) - Q(i-2, i+1) - Q(i-1, i+2)
  \]
  where $Q(j, l) := s^{(n_1, \cdots, n_j)}_{j+1} s^{(n_l, \cdots, n_k)}_{k-l+2} + s^{(n_2, \cdots, n_j)}_{j} s^{(n_l, \cdots, n_{k-1})}_{k-l+1}$.
  By inverting some of $n_*$'s signs with (vii) in Lemma \ref{continuant-lemma1} and applying the Euler identity,
  \begin{align*}
    Q(j, l) & = K(n_1 t, \cdots, n_j t) K(n_l t, \cdots, n_k t) + K(n_2 t, \cdots, n_j t) K(n_l t, \cdots, n_{k-1} t) \\
    & = (-1)^j K(- n_j t, \cdots, - n_1 t) K(n_k t, \cdots, n_l t) \\ & \phantom{K(n_1 t, \cdots, n_j t)} + (-1)^{j-1} K(-n_j t, \cdots, -n_2 t) K(n_{k-1} t, \cdots, n_l t) \\
    & = (-1)^j K(- n_j t, \cdots, - n_1 t, n_k t, \cdots, n_l t).
  \end{align*}
  which gives the desired result for $\hat{T}$.
\end{proof}

\subsection{Roots of ${\bf n}$-Chebyshev polynomials} \label{section:roots}
In this subsection, we investigate properties of the roots of ${\bf n}$-Chebyshev polynomials. Namely, limits of roots of $s^{\bf n}_{k+1}(t)$
and the roots for specific ${\bf n}$. The first property will be described by the limit of $s^{{\bf n}}_{k+1}(t)/n_i$ as a chosen parameter $\abs{n_i}$ tends to infinity.

Note that the coefficient of $n_i$ in (i) of Proposition \ref{split s} is a product of $s^{{\bf n}_L}$ and $s^{{\bf n}_R}$. In general, the following holds.
\begin{lemma} \label{ni limit of s}
    Let ${\bf n} = (n_1, n_2, \cdots, n_k)$ be an integer sequence. For a subsequence of ${\bf n}$, $(n_{i_1}, n_{i_2}, \cdots, n_{i_j})$ with $0 = i_0 < i_1 < i_2 < \cdots < i_m < i_{m+1}=k+1$, 
    $$
        \frac{s^{{\bf n}}_{k+1}(t)}{n_{i_1} n_{i_2} \cdots n_{i_m}} \rightarrow t^m \prod_{j=1}^{m+1} s^{{\bf n}_j}_{{i_{j}} - {i_{j-1}}} (t) \ \ \text{ as } \abs{n_{i_1}}, \abs{n_{i_2}}, \cdots, \abs{n_{i_m}} \rightarrow \infty
    $$ where ${\bf n}_j = (n_{i_{j-1}+1}, \cdots, n_{i_{j}-1})$ for each $j=1,2, \cdots, m+1$.
\end{lemma}
\begin{proof}
    By (i) of Proposition \ref{split s}, $s^{{\bf n}}_{k+1}(t)$ can be written as a linear polynomial with respect to $n_{i_1}$ as follows.
    \begin{align} \label{ni linear decompose}
        s^{{\bf n}}_{k+1}(t) = ts^{{\bf n}_1}(t)s^{(n_{i_1+1}, \cdots, n_k)}(t) \cdot n_{i_1} + g_1(t)
    \end{align}
     where $g_1(t)$ is a polynomial of $t, n_1, \cdots, n_k$ except $n_{i_1}$. Note that $s^{{\bf n}_1}(t)$ does not contain any monomial involving $n_{i_1}, n_{i_2}, \cdots, n_{i_m}$. Therefore, it implies the statement for $m=1$ directly. For $m>1$, one can show that the coefficient of $n_{i_1} n_{i_2} \cdots n_{i_m}$ is $t^m \prod_{j=1}^{m+1} s^{{\bf n}_j}_{{i_{j}} - {i_{j-1}}} (t)
    $ in $s^{{\bf n}}_{k+1}(t)$ and the other terms vanish after dividing by $n_{i_1} n_{i_2} \cdots n_{i_m}$ in the limit. For instance, if $m=2$ then, by applying the same proposition to $s^{(n_{i_1}+1, \cdots, n_k)}(t)$ in \eqref{ni linear decompose}, we have
    $$
        s^{{\bf n}}_{k+1}(t) = t^2s^{{\bf n}_1}(t)s^{{\bf n}_2}(t)s^{{\bf n}_3} \cdot n_{i_1}n_{i_2} + g_1(t) +  tn_{i_1}s^{{\bf n}_1}g_2(t)
    $$ where $g_2(t)$ is a polynomial in $t$ and $n_{i_1+1}, \cdots, n_k$ except $n_{i_1}$ and $n_{i_2}$.
\end{proof}

Now we have the following proposition describing the limiting behaviors of roots of $s^{{\bf n}}_{k+1}(t)$.


\begin{proposition}
    Let $k > 1$ be an integer and $Z_k \subset \mathbb{C}$ be the set of complex numbers $z$ such that $s^{{\bf n}}_{k+1}(z) = 0$ for some nonzero integer sequence ${\bf n}$ of length $k$. Then the derived set of $Z_k$ can be given as
    $$
        Z'_k = \bigcup_{l=1}^{k-1} Z_l.
    $$
\end{proposition}
\begin{proof}
    It suffices to show for nonzero case, i.e., that $Z'_k \setminus \{ 0 \} = \bigcup_{l=1}^{k-1} (Z_l \setminus \{ 0 \})$ holds. 
    Indeed, both sides include $0$: $0 \in Z_1$ by definition, and $0 \in Z'_k$ considering roots of $s_{k+1}^{(c, \cdots, c)}$ converges to $0$ as $c \to \infty$.

    Suppose $z \in Z_l \setminus \{ 0 \}$ for some $l < k$,
    with an integer sequence ${\bf n} = (n_1, n_2, \cdots, n_l)$ such that $s^{{\bf n}}_{l+1}(z) = 0$.
    Then, choose another sequence ${\bf m} = (m_1, \cdots, m_{k-l-1})$ of length $k-l-1$ so that $s^{\bf m}_{k-l}(z) \neq 0$.
    By Lemma \ref{ni limit of s}, $s_{k+1}^{({\bf n}, j, {\bf m})}(t) / j$ converges to $t s_{l+1}^{\bf n}(t) s_{k-l}^{\bf m}(t)$ as $\abs{j} \to \infty$. Thus, by the well-known fact that roots of a polynomial is a continuous function of its coefficients (see \cite{Coolidge} \cite{Gary-Clyde_roots_of_polys} for instance),
    we get a sequence of complex roots $\{ z_j \}_{j=1}^\infty$ converging to $z$, where $z_j \in Z_k$ is the corresponding root of $s_{k+1}^{({\bf n}, j, {\bf m})}(t)$.

    Notice that $z_j \neq z$ is enough to conclude that $z \in Z_k'$, as desired.
    If $z_j = z$ for some $j \geq 1$,
    one has, by Proposition \ref{split s},
    \begin{align*}
        0 &= s_{k+1}^{({\bf n}, j, {\bf m})} (z) \\
        &= j z s_{l+1}^{\bf n}(z) s_{k-l}^{\bf m}(z) - (s_{l+1}^{\bf n}(z) s_{k-l-1}^{(m_2, \cdots, m_{k-l-1})}(z) + s_{l}^{\bf n}(z) s_{k-l}^{\bf m}(z)) \\
        &= - s_l^{\bf n}(z) s_{k-l}^{\bf m}(z)
    \end{align*}
    as $s_{l+1}^{\bf n}(z) = 0$. Since $s_{l+1}^{\bf n}(z)$ and $s_{l}^{\bf n}(z)$ cannot be zero simultaneously by the recursive definition of $s^{{\bf n}}$, this implies $s_{k-l}^{\bf m}(z) = 0$, which contradicts our choice of ${\bf m}$. This proves $Z'_k \setminus \{0\} \supset \bigcup_{l=1}^{k-1} (Z_l \setminus \{0\})$.

       Conversely, suppose $z \in Z_k' \setminus \{ 0 \}$. 
    Take a complex sequence $\{ z_j \}_{j=1}^\infty$ in $Z_k$ which converges to $z$, and consider ${\bf n}^j = (n^j_{1}, \cdots, n^j_{k})$ for each $j$ which is the integer sequence satisfying $s_{k+1}^{{\bf n}^j}(z_j) = 0$.
    Now, consider the sequence $\{ n^j_i \}_{j=1}^\infty$ in terms of $j$ for each $i$.
    Depending on the value of $\limsup_{j \to \infty} \abs{n^j_i}$, either $\abs{n^j_i} \to \infty$ or $\{ n^j_i \}$ is bounded, and we may take a subsequence to assume $n^j_i = n_i$ in the latter case.
    In particular, there is $i_1, \cdots, i_m$ where $\lim_{j \to \infty} \abs{n^j_{i_q}} = \infty$, and $n^j_i = n_i$ is fixed for all the other $i$'s.
    Therefore, we can apply the Lemma \ref{ni limit of s} to see that $s_{k+1}^{{\bf n}^j}(t)$ converges to the product of $s_{i_q - i_{q-1}}^{(n_{i_{q-1}+1}, \cdots, n_{i_q})}(t)$'s and $t^m$.
    Then, $z = \lim_{j \to \infty} z_j$ is a root of the product, by the aforementioned well-known fact of the continuity of roots of polynomials.
    Further, the nonzero $z$ cannot be a root of $t^m = 0$, so $z$ should be a root of $s_{i_q - i_{q-1}}^{(n_{i_{q-1}+1}, \cdots, n_{i_q})}(t) = 0$ for some $1 \leq q \leq m+1$.
    In particular, $z \in Z_l$ for $l = i_q - i_{q-1} < k$. This shows $Z'_k \setminus \{0\} \subset \bigcup_{l=1}^{k-1} (Z_l \setminus \{0\})$, finishing the proof.
\end{proof}

\begin{remark}\label{Smilga}
As stated in Theorem A of the introduction, $\lambda$ is a relation number if and only if 
$\sqrt{-\lambda}$ is a root of an $\bf n$-Chebyshev polynomial $s^{\bf n}_{k+1}(t) \in \mathbb Z[t]$ for some nonzero integer sequence ${\bf n}=(n_1,n_2,\cdots,n_k)$.  This result will be demonstrated in the next section.
    Smilga \cite{Ilia} showed that
    $(\frac{1+\sqrt{5}}{2})^{\pm 2}$ and $2 \pm \sqrt{2}$ are obtained as accumulation points of rational relation numbers. Since the corresponding $\bf n$-sequences have constant length with some $n_i \to \infty$ in each case, we can apply Lemma \ref{ni limit of s} to show that they also are relation numbers.

\end{remark}

Finally, let's consider the roots of $s^{{\bf n}}_{k+1}(t)$ and $T^{{\bf n}}_{k+1}(t)$ for some special sequences ${{\bf n}}$.
If ${\bf n} = (1, 1, 1, \cdots)$, we get the classical Chebyshev polynomials, i.e.
    $$
        s^{\bf n}_{k+1}(t) = U_{k}(\frac{t}{2}) = \prod_{i=1}^k (t - 2\cos \frac{i}{k+1}\pi) \text{  and  } T^{\bf n}_{k}(t) = 2T_{k}(\frac{t}{2}) = \prod_{i=1}^{k} (t - 2\cos \frac{2i-1}{2k}\pi).
    $$
The following is obtained by (i) and (iii) of Lemma \ref{n-cheby-1}, and will be used later.

\begin{lemma} \label{T_alt}
For ${\bf n} = (1,-1, 1, -1,\cdots)$,
 $$
        s^{\bf n}_{2k}(t) = t\prod_{i=1}^{k-1} (-t^2 - 4\cos^2 \frac{i}{2k}\pi) \,,\, T^{\bf n}_{2k}(t)= \prod_{i=1}^{k} (-t^2 - 4\cos^2 \frac{2i-1}{4k}\pi)
    $$ and
    $$
        s^{\bf n}_{2k+1}(t) = \prod_{i=1}^k (-t^2 - 4\cos^2 \frac{i}{2k+1}\pi) \,,\, T^{\bf n}_{2k+1}(t)= t \prod_{i=1}^{k} (-t^2 - 4\cos^2 \frac{2i-1}{4k+2}\pi)
    $$ for $k \geq 0$.
\end{lemma}

\section{Non-freeness of groups generated by two parabolic elements}
\label{section:main-A}
Let $X$ and $Y$ be two parabolic elements in $SL(2,\bc)$, i.e., $tr X=tr Y=2$, such that $XY\neq YX$. Then after conjugating if necessary, we may assume that 
\begin{equation*}
X=\begin{pmatrix}
1& 1\\
0 & 1 
 \end{pmatrix} \quad \text{and}\quad
Y=\begin{pmatrix}
1& 0\\
\lambda & 1 
 \end{pmatrix}=Y_{\lambda}
\end{equation*}
where $\lambda=2-tr(XY^{-1})$, see  \cite{Riley4}. 

\begin{definition}
Let $G_{\lambda}$ be the subgroup of $SL(2,\bc)$ generated by 
$X$ and $Y_{\lambda}$.  The complex number $\lambda$ is called a {\it free number} if $G_{\lambda}$ is a free
group of rank $2$. Otherwise, we say that $\lambda$ is 
 {\it non-free} or 
a {\it relation number}. 
\end{definition}
Note that the terminology `relation number' was introduced by Kim and Koberda \cite{KK}.

Let $F_{x,y} = \langle x,y\rangle$ be a free group of rank two.
If  $\lambda$ is a relation number, then there is a nontrivial word of the form
$$w=x^{n_1}y^{n_2}x^{n_3}y^{n_4}\cdots x^{n_{2m-1}}y^{n_{2m}}\in F_{x,y}$$
for some integer sequence ${\bf n}=(n_1,n_2,\cdots, n_{2m})$
such that $w(X,Y_{\lambda})$ is an identity matrix in $SL(2,\bc)$. (Note that either $n_1$ or $n_{2m}$ is possibly zero.) So we denote $w(X,Y_{\lambda})$ for each integer sequence ${\bf n}$ by $W({\bf n})$ or $W(n_1, n_2,\cdots,n_{2m})$. That is, 
 \begin{equation}
\begin{split}
W({\bf n})=&
X^{n_1}Y_{\lambda}^{n_2}X^{n_3}Y_{\lambda}^{n_4}\cdots X^{n_{2m-1}}Y_{\lambda}^{n_{2m}}\\
=&\begin{pmatrix}
1& 1\\
0 & 1 
 \end{pmatrix}^{n_1}\begin{pmatrix}
1& 0\\
\lambda & 1 
 \end{pmatrix}^{n_2}\begin{pmatrix}
1& 1\\
0 & 1 
 \end{pmatrix}^{n_3}\cdots \begin{pmatrix}
1& 1\\
0 & 1 
 \end{pmatrix}^{n_{2m-1}}\begin{pmatrix}
1& 0\\
\lambda & 1 
 \end{pmatrix}^{n_{2m}}.
\end{split}
\end{equation}

Let $w^*$ be the word  obtained from $w$ by exchanging $x$ with $y$, and $W^*({\bf n}):=w^*(X,Y_{\lambda})$, that is, 
\begin{equation*}
\begin{split}
W^*({\bf n})&=Y_{\lambda}^{n_1}X^{n_2}Y_{\lambda}^{n_3}X^{n_4}\cdots Y_{\lambda}^{n_{2m-1}}X^{n_{2m}}\\
&=\begin{pmatrix}
1& 0\\
\lambda & 1 
 \end{pmatrix}^{n_1}\begin{pmatrix}
1& 1\\
0 & 1 
 \end{pmatrix}^{n_2}\cdots \begin{pmatrix}
1& 0\\
\lambda & 1 
 \end{pmatrix}^{n_{2m-1}}\begin{pmatrix}
1& 1\\
0 & 1 
 \end{pmatrix}^{n_{2m}}.
\end{split}
\end{equation*}
Note that  $CXC^{-1}=Y_{\lambda}$ and $CY_{\lambda}C^{-1}=X$ for $C=\begin{pmatrix}
0& 1\\
\lambda & 0 
 \end{pmatrix}$, and thus  
 \begin{equation}\label{w-w*}
W^*({\bf n})=CW({\bf n})C^{-1}.
\end{equation}
Then we have the following which has already been proved for the case when $n_i \in\{1,-1\}$ in \cite{Jo-Kim2}.
\begin{theorem}\label{W-eq}
Let  ${\bf n}$ be an integer sequence ${\bf n}=(n_1,n_2,\cdots, n_{2m})$ and $\lambda=-u^2$.
We have
\begin{enumerate}
\item[\rm (i)]
$W({\bf n})=(-1)^{m-1}\begin{pmatrix}
-s^{\bf n}_{2m+1}(u) & \displaystyle \frac{1}{u} s^{\bf n}_{2m}(u)\\
-u  s^{\sigma(\bf n)}_{2m}(u)& s^{{\sigma(\bf n)}}_{2m-1}(u)
 \end{pmatrix}$
\item[\rm (ii)]
$ W^*({\bf n})=(-1)^{m-1}\begin{pmatrix} s^{{\sigma(\bf n)}}_{2m-1}(u) & \displaystyle \frac{1}{u}s^{\sigma(\bf n)}_{2m}(u)\\
-u s^{\bf n}_{2m}(u)& 
-s^{\bf n}_{2m+1}(u)\end{pmatrix}$
\end{enumerate}
 \end{theorem}
\begin{proof}
Let $W({\bf n})=\begin{pmatrix}
W_{11}& W_{12}\\
W_{21} & W_{22} 
 \end{pmatrix}$.  
\begin{enumerate}
\item[\rm (i)]
Since
\begin{equation*}
\begin{split}
W({\bf n})
&=\begin{pmatrix}
1& 1\\
0 & 1 
 \end{pmatrix}^{n_1}\begin{pmatrix}
1& 0\\
\lambda & 1 
 \end{pmatrix}^{n_2}\begin{pmatrix}
1& 1\\
0 & 1 
 \end{pmatrix}^{n_3}\cdots \begin{pmatrix}
1& 1\\
0 & 1 
 \end{pmatrix}^{n_{2m-1}}\begin{pmatrix}
1& 0\\
\lambda & 1 
 \end{pmatrix}^{n_{2m}}\\
&=\begin{pmatrix}
1 & n_1\\
  0 & 1
 \end{pmatrix}\begin{pmatrix}
1 & 0\\
 n_2\lambda&1
 \end{pmatrix}\begin{pmatrix}
1 & n_3\\
  0 & 1
 \end{pmatrix}\begin{pmatrix}
1 & 0\\
  n_4\lambda& 1
 \end{pmatrix}
\cdots
\begin{pmatrix}
1 & n_{2m-1}\\
  0 & 1
 \end{pmatrix}\begin{pmatrix}
1 & 0\\
  n_{2m}\lambda& 1
 \end{pmatrix}\\
 &=\begin{pmatrix}
n_1 & 1\\
  1 & 0
 \end{pmatrix}\begin{pmatrix}
n_2\lambda & 1\\
1  & 0
 \end{pmatrix}\begin{pmatrix}
n_3 & 1\\
  1 & 0
 \end{pmatrix}\begin{pmatrix}
n_4\lambda & 1\\
  1& 0
 \end{pmatrix}
\cdots
\begin{pmatrix}
n_{2m-1} & 1\\
  1 & 0
 \end{pmatrix}\begin{pmatrix}
n_{2m}\lambda & 1\\
  1& 0
 \end{pmatrix}\\
&=M^+(n_1,n_2\lambda,n_3,n_4\lambda,\cdots,n_{2m-1},n_{2m}\lambda)\\
&=
\begin{pmatrix}
K^+_{2m}(  n_1, n_{2}\lambda,\cdots,n_{2m-1},n_{2m}\lambda) & K^+_{2m-1}(  n_1, n_{2}\lambda,\cdots,n_{2m-2}\lambda, n_{2m-1}) \\
K^+_{2m-1}(  n_{2}\lambda, n_3\cdots, n_{2m-1}, n_{2m}\lambda) & K^+_{2m-2}(  n_{2}\lambda,n_{3}, \cdots,n_{2m-2}\lambda, n_{2m-1})
\end{pmatrix},\\
\end{split}
\end{equation*}
the statement is proved using  Proposition \ref{Kplus&K} and (i) of  Lemma \ref{s-v-continuant} as follows:
\begin{equation*}
\begin{split}
W_{11}&=K^+_{2m}(  n_1, n_{2}\lambda,n_{3}, n_{4}\lambda,\cdots,n_{2m-2}\lambda, n_{2m-1},n_{2m}\lambda)\qquad\qquad\qquad\\
&=(-1)^mK_{2m}(n_1u,n_2u,n_3u\cdots,n_{2m-1}u,n_{2m}u)=(-1)^ms^{\bf n}_{2m+1}(u)\\
W_{12}&=K^+_{2m-1}(  n_1, n_{2}\lambda,n_{3}, n_{4}\lambda,\cdots,n_{2m-2}\lambda, n_{2m-1})\\
&=(-1)^{m-1}\frac{1}{u}K_{2m-1}(n_1u,n_2u,n_3u\cdots,n_{2m-1}u)=(-1)^{m-1}\frac{1}{u}s^{\bf n}_{2m}(u)\\
W_{21}&=K^+_{2m-1}(  n_{2}\lambda,n_{3}, n_{4}\lambda,\cdots,n_{2m-2}\lambda, n_{2m-1}, n_{2m}\lambda)\\
&=(-1)^muK_{2m-1}(n_2u,n_3u\cdots,n_{2m-1}u,n_{2m}u)=(-1)^mus^{\sigma(\bf n)}_{2m}(u)\\
W_{22}&=K^+_{2m-2}(  n_{2}\lambda,n_{3}, n_{4}\lambda,\cdots,n_{2m-2}\lambda, n_{2m-1})\\
&=(-1)^{m-1}K_{2m-2}(n_2u,n_3u\cdots,n_{2m-1}u)=(-1)^{m-1}s^{\sigma(\bf n)}_{2m-1}(u)\\
\end{split}
\end{equation*}
\item[\rm (ii)]
By (\ref{w-w*}), 
we have 
\begin{equation*}
        W^*({\bf n})=CW({\bf n})C^{-1}=
\begin{pmatrix}
{W}_{22} & \lambda^{-1}W_{21} \\
\lambda W_{12}& {W}_{11}
 \end{pmatrix},\,\,
 C=\begin{pmatrix}
0& 1\\
\lambda & 0 
 \end{pmatrix}.
\end{equation*}
This completes the proof.
\end{enumerate}
\end{proof}

\begin{corollary}
Let  ${\bf n}$ be an integer sequence ${\bf n}=(n_1,n_2,\cdots, n_{2m})$ and $\lambda=-u^2$. Then we have 
\begin{enumerate}
\item[\rm (i)]
$X^{n_1}Y_{\lambda}^{n_2}X^{n_3}\cdots  X^{n_{2m-1}}Y_{\lambda}^{n_{2m}}
=(-1)^{m-1}\begin{pmatrix}
-s^{\bf n}_{2m+1}(u) & \displaystyle \frac{1}{u}s^{\bf n}_{2m}(u)\\
-u  s^{\sigma(\bf n)}_{2m}(u)& s^{\sigma(\bf n)}_{2m-1}(u)
 \end{pmatrix}$
\item[\rm (ii)]
$X^{n_1}Y_{\lambda}^{n_2}X^{n_3}\cdots X^{n_{2m-1}}
=(-1)^{m-1}\begin{pmatrix}
s^{\bf n}_{2m-1}(u) & \displaystyle\frac{1}{u}s^{\bf n}_{2m}(u)\\
u  s^{\sigma(\bf n)}_{2m-2}(u)& s^{\sigma(\bf n)}_{2m-1}(u)
 \end{pmatrix}
$
\item[\rm (iii)]
$
Y_{\lambda}^{n_1}X^{n_2}Y_{\lambda}^{n_3}\cdots Y_{\lambda}^{n_{2m-1}}
=(-1)^{m-1}\begin{pmatrix}
s^{\sigma(\bf n)}_{2m-1}(u) & \displaystyle -\frac{1}{u}s^{\sigma(\bf n)}_{2m-2}(u)\\
-u  s^{\bf n}_{2m}(u)& s^{\bf n}_{2m-1}(u)
 \end{pmatrix}
$
\item[\rm (iv)]
$
Y_{\lambda}^{n_1}X^{n_2}Y_{\lambda}^{n_3}\cdots Y_{\lambda}^{n_{2m-1}}X^{n_{2m}}
=(-1)^{m-1}\begin{pmatrix}
s^{\sigma(\bf n)}_{2m-1}(u) & \displaystyle \frac{1}{u}s^{\sigma(\bf n)}_{2m}(u)\\
-u  s^{\bf n}_{2m}(u)& -s^{\bf n}_{2m+1}(u)
 \end{pmatrix}
$
\end{enumerate}
\end{corollary}
\begin{proof}
The identities (i) and (iv) are in Theorem \ref{W-eq}, and from theses (ii) and (iii) are obtained by letting $n_{2m}=0$.
\end{proof}
\begin{lemma}\label{even}
For any integer sequence ${\bf n}=(n_1,n_2,\cdots,n_{j})$,  $W({\bf n})=Id$ if and only if $W^*({\bf n})=Id$.
\end{lemma}
\begin{proof}
By (\ref{w-w*}), we have $W^*({\bf n})=CW({\bf n})C^{-1}$ for $C=\begin{pmatrix}
0& 1\\
\lambda & 0 
 \end{pmatrix}$. Hence the statement follows.
\end{proof}
\begin{corollary}\label{even-odd}
If a nonzero complex number $\lambda$ is a relation number, then the following hold:
\begin{enumerate}
\item[\rm (i)]There is a nonzero integer sequence ${\bf n}=(n_1,n_2,\cdots,n_{2k})$ such that  $W({\bf n})=X^{n_1}Y_{\lambda}^{n_2}X^{n_3}Y_{\lambda}^{n_4}\cdots Y_{\lambda}^{n_{2k}}=Id$.
\item[\rm (ii)]There is a nonzero integer sequence ${\bf n}=(n_1,n_2,\cdots,n_{2k+1})$ such that  $W({\bf n})=X^{n_1}Y_{\lambda}^{n_2}X^{n_3}Y_{\lambda}^{n_4}\cdots Y_{\lambda}^{n_{2k}}X^{n_{2k+1}}=Id$.
\end{enumerate}
\end{corollary}
\begin{proof}
Due to the non-freeness of $\lambda$, we can choose a a nonzero integer sequence ${\bf n}=(n_1,n_2,\cdots,n_{j})$  such that $W({\bf n})=Id$.
If $j=2i+1$ then 
$$Id=W({\bf n})W^*({\bf n})=W(\tilde{{\bf n}})\,\, \text{for} \,\, \tilde{{\bf n}}= (n_1,n_2,\cdots,n_{2i+1},n_1,n_2,\cdots,n_{2i+1}),$$
and if $j=2i$ then
$$Id=W({\bf n})W^*({\bf n})=W(\bar{{\bf n}})\,\, \text{for} \,\, \bar{{\bf n}}= (n_1,n_2,\cdots,n_{2i}+n_1,n_2,\cdots,n_{2i}).$$
Note that the sequence $\tilde{{\bf n}}$ is a nonzero integer sequence of even length, and  the sequence $\bar{{\bf n}}$ might be able to be reduced  but the length of  $\bar{{\bf n}}$ is odd whether or not reduced.
This completes the proof.
\end{proof}

We are now in a position to derive a necessary and sufficient condition for a complex number $\lambda$ to be a relation number, as stated in the introduction in Theorem \ref{main2}. First we can prove the following results regarding the ${\bf n}$-Chebyshev polynomial $s^{\bf n}_k(t)$, a generalized Chebyshev polynomial of the second kind. 
\begin{theorem}\label{rel number}
Let  $\lambda$ be a nonzero complex number. Then the followings are all equivalent:
\begin{enumerate}
\item[\rm (i)]
 $\lambda$ is a relation number.
\item[\rm (ii)] 
$\sqrt{-\lambda}$ is a root of an $\bf n$-Chebyshev polynomial $s^{\bf n}_{k+1}(t) \in \mathbb Z[t]$ for some nonzero integer sequence ${\bf n}=(n_1,n_2,\cdots,n_k)$.
\item[\rm (iii)]
$\sqrt{-\lambda}$ is a root of an $\bf n$-Chebyshev polynomial $s^{\bf n}_{2m}(t) \in \mathbb Z[t]$ for some nonzero integer sequence ${\bf n}=(n_1,n_2,\cdots,n_{2m-1})$.
\end{enumerate}
\end{theorem}
\begin{proof}
We first prove that (ii) and (iii) are equivalent. To do this, we just need to prove the following:
$$s^{\bf n}_{k}(u)=0  \,\, \text{implies}\,\,\, s^{\bar{\bf  n}}_{2k}(u)=0\,\, \text{ for some  integer sequence}\,\,\bar{\bf  n}.$$
Let ${\bf  n}=(n_1, \cdots, n_{k-1})$. 
If we choose  $\bar{\bf  n}$  as
$\bar{\bf  n}=(n_1, \cdots, n_{k-1}, 1, n_1, \cdots, n_{k-1})$, then $s^{\bar{\bf  n}}_{2k}(u)=0$
by Proposition \ref{EI-1}. Actually   Proposition \ref{EI-1} implies that any of the following can be $\bar{\bf  n}$ for any nonzero integer $x$:
$$({\bf  n}, x, {\bf  n}),\,({\bf  n}, x, -{\bf  n}),\,({\bf  n},x, \overleftarrow{{\bf n}}), ({\bf  n},x, -\overleftarrow{{\bf n}}),$$
where $\overleftarrow{{\bf n}}=(n_{k-1},n_{k-2},\cdots,n_{2},n_1)$. This proves that  (ii) and (iii) are equivalent. So it suffices to show that a complex number  $\lambda$ is a relation number if and only if (iii) holds.

If  $\lambda$ is a relation number, then there exists a nonzero integer sequence ${\bf n}=(n_1,n_2,\cdots, n_{2m})$ such that $W({\bf n})=X^{n_1}Y_{\lambda}^{n_2}X^{n_3}Y_{\lambda}^{n_4}\cdots X^{n_{2m-1}}Y_{\lambda}^{n_{2m}}=Id$ by Corollary \ref{even-odd}. This implies $s^{\bf n}_{2m}(\sqrt{-\lambda})=0$, since
$$W({\bf n})=(-1)^{m-1}\begin{pmatrix}
-s^{\bf n}_{2m+1}(\sqrt{-\lambda}) & \displaystyle \frac{s^{\bf n}_{2m}(\sqrt{-\lambda})}{\sqrt{-\lambda}}\\
- \sqrt{-\lambda} \,\, s^{\sigma(\bf n)}_{2m}(\sqrt{-\lambda})& s^{\sigma(\bf n)}_{2m-1}(\sqrt{-\lambda})
 \end{pmatrix}$$
 by Theorem \ref{W-eq}.
Conversely, if $s^{\bf n}_{2m}(\sqrt{-\lambda})=0$, then 
$$W({\bf n})=\begin{pmatrix}
* & 0\\
*  & *
 \end{pmatrix}$$
by Theorem \ref{W-eq} again, and thus $[WY_{\lambda}W^{-1},Y_{\lambda}]=Id$ by the Lemma \ref{LU} below, which implies that $G_{\lambda}$ is not a free group of rank $2$ and thus $\lambda$ is a relation number.
\end{proof}
\begin{lemma}[\cite{LU}]\label{LU}
 Let $A$ be a lower triangular matrix in $SL(2,\bc)$. Then $AY_{\lambda}A^{-1}$ commutes with $Y_{\lambda}$ for any $\lambda\in \bc$.
\end{lemma}
\begin{proof}
Let $A=\begin{pmatrix}
a& 0\\
b & c 
 \end{pmatrix}\in SL(2,\bc)$. Then $[AY_{\lambda}A^{-1},Y_{\lambda}]=Id$ follows from 
$$AY_{\lambda}A^{-1}=
\begin{pmatrix}
a& 0\\
b & c 
 \end{pmatrix}
 \begin{pmatrix}
1& 0\\
\lambda & 1 
 \end{pmatrix}
 \begin{pmatrix} c& 0\\   -b & a  \end{pmatrix}
=\begin{pmatrix}
1& 0\\
c^2\lambda & 1 
 \end{pmatrix}.$$
\end{proof}
 From the identity $s^{(\overleftarrow{{\bf n}},{\bf n})}_{2k+1}(t)=(-1)^kv^{\bf n}_{k}(t)v^{\bf n}_{k}(-t)$ for any integer sequence ${\bf n}=(n_1,n_2,\cdots,n_k)$ in Proposition \ref{factor-cont-prop}, we immediately  obtain the following:
\begin{corollary}\label{3rd-4th}
A nonzero complex number  $\lambda$ is a relation number if  $\sqrt{-\lambda}$ is a root of either $v^{\bf n}_{k}(t)$ or  $v^{\bf n}_{k}(-t)$ for some integer sequence ${\bf n}$ and some integer $k>0$.
\end{corollary}
Recall that $v^{\bf n}_{k}(t)$ and  $v^{\bf n}_{k}(-t)$ are generalized Chebyshev polynomials of the third and fourth kind, respectively.  
Now we show that all the roots of $T^{\bf n}_k, \tilde{T}^{\bf n}_k$, and $\hat{T}^{\bf n}_k$, generalized Chebyshev polynomials of the first kind, are also relation numbers, and vice versa.

\begin{lemma} \label{T-lemma}
If a nonzero complex number $\lambda$ is a relation number, then the following holds:
\begin{enumerate}
\item[\rm (i)]
There is a nonzero integer sequence ${\bf n}=(n_1,n_2,\cdots,n_{k})$ such that  
$$\tilde{T}^{\bf n}_{k}(\sqrt{-\lambda})=0=\hat{T}^{\bf n}_{k}(\sqrt{-\lambda})$$
\item[\rm(ii)]There is a nonzero integer sequence ${\bf n}=(n_1,n_2,\cdots,n_{k})$ such that  
$$\tilde{T}^{\bf n}_{k}(\sqrt{-\lambda})=\pm 2 \,\,\text{and}\,\, \hat{T}^{\bf n}_{k}(\sqrt{-\lambda})=0$$
\end{enumerate}
\end{lemma}
\begin{proof}
Let $\lambda=-u^2$. Since $\lambda$ is a relation number,  there are two integer sequences ${\bf n}=(n_1,n_2,\cdots,n_{2i+1})$ and ${\bf m}=(m_1,m_2,\cdots,m_{2i})$ such that $W({\bf n})=Id=W({\bf m})$ by Corollary \ref{even-odd}.
For the sequence $\bf n$, we have 
$$Id=W({\bf n})
=(-1)^{i}\begin{pmatrix}
s^{\bf n}_{2i+1}(u) & \displaystyle\frac{1}{u}s^{\bf n}_{2i+2}(u)\\
u  s^{\sigma(\bf n)}_{2i}(u)& s^{\sigma(\bf n)}_{2i+1}(u)
 \end{pmatrix},
$$
and this implies that
$$\tilde{T}^{\bf n}_{2i+1}(u)=s^{\bf n}_{2i+2}(u)- s^{\sigma(\bf n)}_{2i}(u)=0-0=0$$
and
$$\hat{T}^{\bf n}_{2i+1}(u)=s^{\bf n}_{2i+2}(u)+ s^{\sigma(\bf n)}_{2i}(u)=0+0=0.$$
 For the sequence $\bf m$, we have 
$$Id=W({\bf m})
=(-1)^{i-1}\begin{pmatrix}
-s^{\bf m}_{2i+1}(u) & \displaystyle \frac{1}{u}s^{\bf m}_{2i}(u)\\
-u  s^{\sigma(\bf m)}_{2i}(u)& s^{\sigma(\bf m)}_{2i-1}(u)
 \end{pmatrix},$$
and this implies that
$$\tilde{T}^{\bf m}_{2i}(u)=s^{\bf m}_{2i+1}(u)- s^{\sigma(\bf m)}_{2i-1}(u)=(-1)^{i}-(-1)^{i-1}=(-1)^i 2$$
and
$$\hat{T}^{\bf m}_{2i}(u)=s^{\bf m}_{2i+1}(u)+ s^{\sigma(\bf m)}_{2i-1}(u)=(-1)^{i}+(-1)^{i-1}=0.$$

\end{proof}

\begin{lemma}\label{lem-Bam}
For any integer sequence ${\bf n}=(n_1,n_2,\cdots,n_{2k+1})$, the following hold:
\begin{enumerate}
\item[\rm (i)]
 $W^*({\bf n})=W({\bf n})^{-1}\,\,\text{ if and only if}\,\, \tilde{T}^{\bf n}_{2k+1}(u)=0$
\item[\rm (ii)]
 $W^*({\bf n})=W({\bf -n})^{-1}\,\,\text{ if and only if}\,\, \hat{T}^{\bf n}_{2k+1}(u)=0$
\end{enumerate}
On the other hand, for ${\bf n}=(n_1,n_2,\cdots,n_{2k})$, we have
$$W({\bf n})=W({\bf -n})^{-1}\text{ if and only if}\,\, \hat{T}^{\bf n}_{2k}(u)=0.$$
\end{lemma}
\begin{proof}
From
$$
W^*({\bf n})=(-1)^{k}\begin{pmatrix}
s^{\sigma(\bf n)}_{2k+1}(u) & \displaystyle -\frac{1}{u}s^{\sigma(\bf n)}_{2k}(u)\\
-u  s^{\bf n}_{2k+2}(u)& s^{\bf n}_{2k+1}(u)
 \end{pmatrix}
$$
and 
$$W({\bf n})^{-1}=(-1)^{k}\begin{pmatrix}
s^{\sigma(\bf n)}_{2k+1}(u) & \displaystyle -\frac{1}{u}s^{\bf n}_{2k+2}(u)\\
-u  s^{\sigma(\bf n)}_{2k}(u)& s^{\bf n}_{2k+1}(u)
 \end{pmatrix}$$
it is obvious that 
$W^*({\bf n})=W({\bf n})^{-1}\,\,\text{ if and only if}\,\, \tilde{T}^{\bf n}_{2k+1}(u)=s^{\bf n}_{2k+2}(u)- s^{\sigma(\bf n)}_{2k}(u)=0$. This proves (i). (ii) and the last statement can also be similarly proven.
\end{proof}

\begin{proposition}\label{TandS}
A nonzero complex number  $\lambda$ is a relation number if  $\sqrt{-\lambda}$ is a root of any one of three polynomials: 
$T^{\bf n}_{k}(t)$, $\tilde{T}^{\bf n}_{k}(t)$, and $\hat{T}^{\bf n}_{k}(t)$, where
 ${\bf n}=(n_1,n_2,\cdots, n_{k})$ is a nonzero  integer sequence.
\end{proposition}
\begin{proof}
Let ${\bf m}=(n_1,\cdots, n_{k},n_{k-1},\cdots,n_1)$, $\tilde{{\bf m}}=(n_1,\cdots, n_{k},n_{1},\cdots,n_{k-1})$, and $\hat{{\bf m}}=(n_1,\cdots, n_{k},-n_1,\cdots, -n_{k-1})$. Then by Proposition \ref{factor-cont-prop},
 $T^{\bf n}_{k}(t)\,|\,s^{{\bf m}}_{2k}(t)$, 
$\tilde{T}^{\bf n}_{k}(t)\,|\,s^{\tilde{{\bf m}}}_{2k}(t)$, and
$\hat{T}^{\bf n}_{k}(t)\,|\,s^{\hat{{\bf m}}}_{2k}(t)$.
Hence the statement follows from Theorem \ref{rel number}.
\end{proof}

From Lemma \ref{T-s identity},  Lemma \ref{T-lemma}, and Proposition \ref{TandS},  we have
\begin{theorem}\label{factor-cont}
Let $\lambda$ be a nonzero complex number. Then 
\begin{enumerate}
\item[\rm (i)]
$\lambda$ is a relation number if and only if  $\sqrt{-\lambda}$ is a root of $T^{\bf n}_{k}(t)$
 for a nonzero  integer sequence ${\bf n}=(n_1,n_2,\cdots, n_{k})$.
\item[\rm (ii)]
$\lambda$ is a relation number if and only if  $\sqrt{-\lambda}$ is a root of $\tilde{T}^{\bf n}_{k}(t)$
 for a nonzero  integer sequence ${\bf n}=(n_1,n_2,\cdots, n_{k})$.
\item[\rm (iii)]
$\lambda$ is a relation number if and only if  $\sqrt{-\lambda}$ is a root of $\hat{T}^{\bf n}_{k}(t)$
 for a nonzero  integer sequence ${\bf n}=(n_1,n_2,\cdots, n_{k})$.
\end{enumerate}
\end{theorem}
\begin{remark}\label{Jang-Kim}
Note that $\tilde{T}^{\bf n}_{2k+1}(u)$ corresponds to Bamberg polynomial in \cite{Bamberg}, since $$W^*({\bf n})=CW({\bf n})C^{-1}$$
and 
$$W^*({\bf n})=W({\bf n})^{-1}\,\,\text{ if and only if}\,\, \tilde{T}^{\bf n}_{2k+1}(u)=0$$
by Lemma \ref{lem-Bam}. Similarly,  one can show that 
$\hat{T}^{\bf n}_{2k+1}(u)$ corresponds to Smilga polynomial in \cite{Ilia}.

On the other hand, $T^{\bf n}_{2k+1}(t)$ is related to  Jang-Kim polynomial discussed  in \cite{Jang-Kim}.  To be more precise, the polynomial $p_k(\alpha)$ constructed in the proof of their main theorem can be associated with our polynomial as follows:
$$p_k(-t^2) =\frac{1}{t}T^{\bf n}_{2k+1}(t), \,\,{\bf n}=(1,-1,1,-1,\cdots).$$
According to Lemma \ref{T_alt}, we can express $p_k(\alpha)$ as:
$$
    p_k(\alpha) = \prod_{j=1}^k (\alpha - 4\cos^2 \frac{2j-1}{4k + 2}\pi)
    $$
The roots of $p_k(\alpha)$, specifically $4\cos^2 \frac{2j-1}{4k + 2}\pi$, are all distinct relation numbers. Moreover, as pointed in \cite{Jang-Kim}, these roots are accumulating relation numbers converging to $4$ as $k \rightarrow \infty.$ 
\end{remark}
Expanding upon Proposition \ref{TandS}, the following much more complicated polynomials also produce relation numbers.
\begin{theorem} \label{continuant_Ts}
Let  ${\bf n}=(n_1,n_2,\cdots,n_k)$ and ${\bf m}=(m_1,m_2,\cdots)$ be nonzero integer sequences.
A nonzero complex number 
 $\lambda$ is a relation number
  if $\sqrt{-\lambda}$ is a root of any one of the following polynomials for $j=1,2,\cdots$:
\begin{enumerate}
\item[\rm (i)]$K_j\left(\tilde{T}^{({\bf n},m_1)}_{k+1}(t), \cdots,\tilde{T}^{({\bf n},m_j)}_{k+1}(t)\right)$
\item[\rm (ii)] $K_{2j}\left(T^{({\bf n},m_1)}_{k+1}(t),T^{(\overleftarrow{{\bf n}},m_2)}_{k+1}(t), \cdots,T^{(\overleftarrow{{\bf n}},m_{2j})}_{k+1}(t)\right)$
\item[\rm (iii)] $K_{2j+1}\left(T^{({\bf n},m_1)}_{k+1}(t),T^{(\overleftarrow{{\bf n}},m_2)}_{k+1}(t), \cdots,T^{({\bf n},m_{2j+1})}_{k+1}(t)\right)$
\item[\rm (iv)] $K_{2j}\left(\hat{T}^{({\bf n},m_1)}_{k+1}(t),\hat{T}^{(\overleftarrow{{\bf n}},m_2)}_{k+1}(t), \cdots,\hat{T}^{(\overleftarrow{{\bf n}},m_{2j})}_{k+1}(t)\right)$
\item[\rm (v)] $K_{2j+1}\left(\hat{T}^{({\bf n},m_1)}_{k+1}(t),\hat{T}^{(\overleftarrow{{\bf n}},m_2)}_{k+1}(t), \cdots,\hat{T}^{({\bf n},m_{2j+1})}_{k+1}(t)\right)$.
\end{enumerate} 
Note that $$\tilde{T}^{({\bf n},m_i)}_{k+1}(t)=m_i t s^{\bf n}_{k+1}(t)-s^{\bf n}_{k}(t)-s^{\sigma({\bf n})}_{k}(t)=s^{({\bf n},m_i)}_{k+2}(t) -s^{\sigma({\bf n})}_{k}(t)$$
$$\hat{T}^{({\bf n},m_i)}_{k+1}(t)=m_i t s^{\bf n}_{k+1}(t)-s^{\bf n}_{k}(t)+s^{\sigma({\bf n})}_{k}(t)=s^{({\bf n},m_i)}_{k+2}(t)+s^{\sigma({\bf n})}_{k}(t)
$$
$$T^{({\bf n},m_i)}_{k+1}(t)=m_i t s^{\bf n}_{k+1}(t)-2s^{\bf n}_{k}(t)=s^{({\bf n},m_i)}_{k+2}(t) -s^{\bf n}_{k}(t).$$
\end{theorem}
\begin{proof}
We get the statement by applying Proposition \ref{cont-fact} to 
$c=(n_1t,\cdots, n_kt)$ 
and
$(x_1, x_2, \cdots)=(m_1t,m_2t,\cdots).$
\end{proof}

\begin{remark}
    One can quickly check that $\sqrt{2}$,   $\sqrt{-2}$, and $\sqrt{3}$ are relation numbers  by observing that $$
        \tilde{T}^{(1,-1,-1,-1)}_4(t) = -t^4 + 2,\, \tilde{T}^{(1,1,-1,-1)}_4(t) = t^4 + 2
    $$
    and
    $$
        K_{2}\left(T^{(1,1,1)}_{3}(t),T^{(1,1,-1)}_{3}(t)\right) = -t^4 + 3
    $$
\end{remark}

\section{Rational relation numbers}
\label{section:rational}
It has been conjectured  that every rational number $q\in(-4,4)$ is a relation number.
 In this section, we present some explicit results in relation to this conjecture as consequences of Theorem \ref{rel number}, (ii).
 \begin{definition}
 A nonzero complex number $\lambda$ is called an $l$-{\it step relation number}
 if $s^{\bf n}_{2l+2}(\sqrt{-\lambda})=0$ for some integer sequence $\bf n$.
 \end{definition}
This is equivalent to that 
 $s^{\bf n}_{2m+2}(\sqrt{-\lambda})=0$ for some $m\in [0,l]$ and some integer sequence $\bf n$,
since  $s^{\bf n}_{2k}(\sqrt{-\lambda})=0$ implies that  $s^{\bf n'}_{2k+2}(\sqrt{-\lambda})=0$ for any integer sequence ${\bf n'}$ such that
$$n'_i=n_i \,\,(i\leq 2k-1),\,\, n'_{2k+1}=0.$$ 
 Note that this definition is consistent with that of Kim and Koberda \cite{KK}. 
 In fact, they defined $l$-step relation number as follows:
 A complex number $\lambda$ is called an $l$-step relation number  if there exists a nontrivial word of the form 
$
w(x,y)=x^{n_1}y^{n_2}x^{n_3}y^{n_4}\cdots x^{n_{2k+1}} 
$
in a free group $F=\langle x,y\rangle$, for some $k\in[0,l]$, such that 
$$w(X,Y_{\lambda})=\begin{pmatrix}
1& 1\\
0 & 1 
 \end{pmatrix}^{n_1}\begin{pmatrix}
1& 0\\
\lambda & 1 
 \end{pmatrix}^{n_2}\begin{pmatrix}
1& 1\\
0 & 1 
 \end{pmatrix}^{n_3}\cdots \begin{pmatrix}
1& 1\\
0 & 1 
 \end{pmatrix}^{n_{2k}}\begin{pmatrix}
1& 0\\
\lambda & 1 
 \end{pmatrix}^{n_{2k+1}}$$ is a lower-triangular matrix in $SL(2,\bc)$. 
 Note that every non-free number is an $l$-step relation number for some $l\geq 0$,  and conversely, according to Lemma \ref{LU}, every $l$-step relation number is a non-free number having a word $w(x,y)\in F$ of syllable length at most $8(l+1)$ such that $w(X,Y_{\lambda})=1$.

 Note that $-\sqrt{-\lambda}$ is also a root of $s^{\bf n}_{2m}(t)$ if $\sqrt{-\lambda}$ is, since $t| s^{\bf n}_{2m}(t)$ and $\frac{1}{t}s^{\bf n}_{2m}(t)$ is an even polynomial. The following observations are well-known.
\begin{proposition}\label{property1}
The following hold.
\begin{enumerate}
\item[\rm (i)] If  a complex number $\lambda$ is an $l$-step relation number, so is $-\lambda$.
\item[\rm (ii)] If $\lambda$ is an $l$-step relation number, so is $\frac{\lambda}{r}$ for any nonzero integer $r$.
\end{enumerate}
\end{proposition}
\begin{proof}
\begin{enumerate}
\item[\rm (i)] The statement follows from 
$\begin{pmatrix}
1& 0\\
-\lambda & 1
 \end{pmatrix}^{-n}=\begin{pmatrix}
1& 0\\
\lambda & 1
 \end{pmatrix}^{n}$.
\item[\rm (ii)]  The statement follows from 
$\begin{pmatrix}
1& 0\\
\frac{\lambda}{r} & 1
 \end{pmatrix}^{rn}=\begin{pmatrix}
1& 0\\
\lambda & 1
 \end{pmatrix}^{n}$.
\end{enumerate}
\end{proof}
Let us denote the set of  complex numbers which  are $l$-step relation numbers  by $R^{(l)}$, that is, 
$$R^{(l)}:=\{q\in \mathbb C \,|\, q\,\, \text{is an $l$-step relation number}\}=\{q\in \mathbb C \,|\, s^{\bf n}_{2l+2}(\sqrt{-q})=0\,\, \text{for some }\,{\bf n}\},$$
and denote the set of  rational numbers which  are $l$-step relation numbers  by $R^{(l)}_{\mathbb Q}$, that is, 
$$R^{(l)}_{\mathbb Q}:=R^{(l)}\cap \mathbb Q. $$
The following proposition characterizes the set of $1$-step relation numbers.
\begin{proposition}\label{property-1step}
Let $A^{(1)}_{\mathbb Q}:=\{\frac{a+b}{ab} \mid a,b,c \in\mathbb Z\setminus\{0\}\}$. Then we have the following:
$$R^{(1)}=R^{(1)}_{\mathbb Q}=\{\frac{a+b}{abc} \,|\, a,b,c \in\mathbb Z\setminus\{0\}\}=A^{(1)}_{\mathbb Q}$$
Consequently, $\frac{1}{r}$ and $ \frac{2}{r}$ are $1$-step relation numbers for any nonzero integer $r$.
\end{proposition}
\begin{proof}
For any nonzero integer sequence ${\bf n}=(n_1,n_2,n_3)$, 
$$
s^{\bf n}_4(u)=u\left (n_1n_2n_3u^2-(n_1+n_3)\right )
$$
and $u$ is a nonzero root of  $s^{\bf n}_4(u)=0$ if and only if $-u^2=-\frac{n_1+n_3}{n_1n_2n_3}$, which implies that $\frac{n_1+n_3}{n_1n_2n_3}$ is an $1$-step relation number for any integer sequence ${\bf n}=(n_1, n_2, n_3, \cdots)$, and such rational numbers of this form cover all  $1$-step relation numbers, that is, 
$$R^{(1)}_{\mathbb Q}=\{\frac{a+b}{abc} \,|\, a,b,c \in\mathbb Z\setminus\{0\}\}=R^{(1)}.$$

 By choosing $c=1$, we obtain 
$$A^{(1)}_{\mathbb Q}=\{\frac{a+b}{ab} \,|\, a,b,c \in\mathbb Z\setminus\{0\}\}\subset R^{(1)}_{\mathbb Q}$$
and  the identity
$$ \frac{a+b}{abc}=\frac{ac+bc}{(ac)(bc)} $$  implies 
$R^{(1)}_{\mathbb Q}\subset A^{(1)}_{\mathbb Q}$. Therefore we have  $R^{(1)}=A^{(1)}_{\mathbb Q}$ as desired. 
Since $1$ is obtained by choosing $a=b=2$, and $2$  is obtained by choosing $a=b=1$, 
$$1, 2 \in A^{(1)}_{\mathbb Q}$$ and thus 
the last statement is proved.
\end{proof}


For $2$-step relation numbers, we obtain some families of rational relation numbers by factorizing $s^{\bf n}_6(u)$.
\begin{proposition}\label{property-2step}
For arbitrary nonzero integers $r,s,t,v, w$, the following holds:
\begin{enumerate}
\item[\rm (i)] 
 $ \frac{3}{r},\, \frac{r+s+t}{rst},\, \frac{r^2s+tv^2}{r^2stv+rstv^2}  \in R^{(2)}_{\mathbb Q}$.
\item[\rm (ii)] 
$\frac{r+w+t}{rst} \in R^{(2)}_{\mathbb Q}$ if $w|rs$.
\item[\rm (iii)] 
$\frac{r+t+v+w}{rst}\in R^{(2)}_{\mathbb Q}$ if $vw\,|\,rst$ and $rv+tv+tw=0$.
\end{enumerate}
\end{proposition}
\begin{proof}
For any nonzero integer sequence ${\bf n}=(n_1,n_2,n_3, n_4,n_5)$, 
$$
s^{\bf n}_6(u)=u\left (n_1n_2n_3n_4n_5u^4-(n_1n_2n_3+n_1n_2n_5+n_1n_4n_5+n_3n_4n_5)u^2+(n_1+n_3+n_5)\right ).
$$
and thus  $u$ is a nonzero root of 
$s^{\bf n}_6(u)=0$ if and only if $\lambda=-u^2$ satisfies $$n_1n_2n_3n_4n_5\lambda^2+(n_1n_2n_3+n_1n_2n_5+n_1n_4n_5+n_3n_4n_5)\lambda+(n_1+n_3+n_5)=0.$$  
In the case when
$n_1n_2=n_4n_5$, we have  
\begin{equation*}
\begin{split}
n_1&n_2n_3n_4n_5\lambda^2+(n_1n_2n_3+n_1n_2n_5+n_1n_4n_5+n_3n_4n_5)\lambda+(n_1+n_3+n_5)\\
 &=(n_1n_2\lambda+1)(n_1n_2n_3\lambda+n_1+n_3+n_5).
\end{split}
\end{equation*}
So we have  $\frac{n_1+n_3+n_5}{n_1n_2n_3}\in R^{(2)}_{\mathbb Q}$ under the condition of $n_5|n_1n_2$. By choosing $n_1=r, n_2=s, n_3=t, n_5=w$, we have the statement (ii),
$$\frac{r+w+t}{rst} \in R^{(2)}_{\mathbb Q}\,\, \text{if}\,\,w|rs.$$
By choosing $n_1=r, n_2=s=n_5, n_3=t$, we have $\frac{r+s+t}{rst} \in R^{(2)}_{\mathbb Q}$. Since $\frac{r+s+t}{rst}=3$ for $r=s=t=1$, we also have $3\in R^{(2)}_{\mathbb Q}$ and thus $\frac{3}{r}\in R^{(2)}_{\mathbb Q}$ by (ii) of Proposition \ref{property1}.

Now we claim  $\frac{r^2s+tv^2}{r^2stv+rstv^2}  \in R^{(2)}_{\mathbb Q}$, which completes the proof of (i).
 If we assume that $n_1+n_3+n_5=0$, then we have 
\begin{equation*}
\begin{split}
n_1&n_2n_3n_4n_5\lambda^2+(n_1n_2n_3+n_1n_2n_5+n_1n_4n_5+n_3n_4n_5)\lambda+(n_1+n_3+n_5)\\
&=n_1n_2n_3n_4n_5\lambda\left (\lambda+\frac{n_1n_2n_3+n_1n_2n_5+n_1n_4n_5+n_3n_4n_5}{n_1n_2n_3n_4n_5}\right )\\
&=n_1n_2n_3n_4n_5\lambda\left (\lambda+\frac{n_1^2n_2+n_4n_5^2}{n_1n_2(n_1+n_5)n_4n_5}\right ),\\
\end{split}
\end{equation*}
Since we obtain $\frac{n_1^2n_2+n_4n_5^2}{n_1n_2(n_1+n_5)n_4n_5}=\frac{r^2s+tv^2}{r^2stv+rstv^2} $ by choosing $n_1=r, n_2=s, n_3=-r-v, n_4=t, n_5=v$, the claim is proved.

By choosing $n_1=r, n_2=s, n_3=t+v, n_4=\frac{rst}{vw}, n_5=w$, (iii) is proved as follows:
\begin{equation*}
\begin{split}
n_1&n_2n_3n_4n_5\lambda^2+(n_1n_2n_3+n_1n_2n_5+n_1n_4n_5+n_3n_4n_5)\lambda+(n_1+n_3+n_5)\\
&=rst\left (rs+\frac{rst}{v}\right )\lambda^2+\left (rs(t+v+w)+(r+t+v)\frac{rst}{v}\right )\lambda+r+t+v+w\\
&=\left (rst\lambda+(r+t+v+w)\right )\left ((rs+\frac{rst}{v})\lambda+1\right )\\
\end{split}
\end{equation*}
\end{proof}

By applying (ii) of Proposition \ref{property-2step} to $w=rs$ and $t=1$, we have that  $1+\frac{1}{s}+\frac{1}{rs}$ is a relation number, and thus for any nonzero integer $v$,
$\frac{1}{v}(1+\frac{1}{s}+\frac{1}{rs})=\frac{1}{v}+\frac{1}{vs}+\frac{1}{vrs}$ is also a relation number
   by (ii) of Proposition \ref{property1}
Hence we get the following.
\begin{corollary}
For nonzero integers $r,s,t$,  the following holds.
$$\frac{1}{r}+\frac{1}{rs}+\frac{1}{rst}\in R^{(2)}_{\mathbb Q}$$
\end{corollary}
For $3$-step relation numbers, we also obtain a family of rational relation numbers by factorizing $s^{\bf n}_8(u)$.
\begin{proposition}\label{property-3step}
For arbitrary nonzero integers $r,s,t,v, w$, we have 
 $$\frac{r+t+v+w}{rst}\in R^{(3)}_{\mathbb Q}$$
if $rst=vw$ and $r^2s+w^2+rst=(r+t)(v+w)$.
\end{proposition}
\begin{proof}
For any nonzero integer sequence ${\bf n}=(n_1,n_2,n_3, n_4,n_5,n_6,n_7)$, 
$$
s^{\bf n}_8(u)=u (u^6n_1n_2n_3n_4n_5n_6n_7-u^4\sum_{\shortstack{$\scriptstyle i_1:\text{odd} $\\$\scriptstyle  i_5\leq 7$}}^{\wedge} n_{i_1}\cdots n_{i_5}+u^2\sum_{\shortstack{$\scriptstyle i_1:\text{odd} $\\$\scriptstyle  i_3\leq 7$}}^{\wedge} n_{i_1}n_{i_2}n_{i_3}-(n_1+n_3+n_5+n_7) ),
$$
and thus  $u$ is a nonzero root of 
$s^{\bf n}_8(u)=0$ if and only if $\lambda=-u^2$ satisfies
$$\lambda^3n_1n_2n_3n_4n_5n_6n_7+\lambda^2\sum_{\shortstack{$\scriptstyle i_1:\text{odd} $\\$\scriptstyle  i_5\leq 7$}}^{\wedge}  n_{i_1}\cdots n_{i_5}+\lambda\sum_{\shortstack{$\scriptstyle i_1:\text{odd} $\\$\scriptstyle  i_3\leq 7$}}^{\wedge} n_{i_1}n_{i_2}n_{i_3}+(n_1+n_3+n_5+n_7)=0.$$
Here $\displaystyle\sum_{\shortstack{$\scriptstyle i_1:\text{odd} $\\$\scriptstyle  i_k\leq n$}}^{\wedge} $ means the following:
\begin{equation*}
\displaystyle\sum_{\shortstack{$\scriptstyle i_1:\text{odd} $\\$\scriptstyle  i_k\leq m$}}^{\wedge}  :0< i_1< i_2 < i_3 <\cdots<i_k\leq m, \,\, i_{2i+1} \text{ is odd},\,\,  i_{2i} \text{ is even}\\
\end{equation*}
In fact, $\displaystyle\sum_{\shortstack{$\scriptstyle i_1:\text{odd} $\\$\scriptstyle  i_k\leq m$}}^{\wedge} n_{i_1}\cdots n_{i_k}$ is the sum of all possible products  obtained by removing $\frac{m-k}{2}$ number of the adjacent pairs $n_j n_{j+1}$ from the full product $n_1 n_2\cdots n_m$.

In the case when 
\begin{equation}\label{3-step-1}
n_1n_2n_3=n_5n_6n_7\,\,\text{and}\,\,
n_1^2n_2+n_6n_7^2+n_1n_2n_3=(n_1+n_3)n_4(n_5+n_7),
\end{equation}
we have
\begin{equation*}
\begin{split}
\lambda^3&n_1n_2n_3n_4n_5n_6n_7+\lambda^2\sum_{\shortstack{$\scriptstyle i_1:\text{odd} $\\$\scriptstyle  i_5\leq 7$}}^{\wedge}  n_{i_1}\cdots n_{i_5}+\lambda\sum_{\shortstack{$\scriptstyle i_1:\text{odd} $\\$\scriptstyle  i_3\leq 7$}}^{\wedge} n_{i_1}n_{i_2}n_{i_3}+(n_1+n_3+n_5+n_7)\\
&=(n_1n_2n_3\lambda +n_1+n_3+n_5+n_7)\left (n_4n_5n_6n_7\lambda^2+(n_1n_2+n_6n_7)\lambda+1\right )
\end{split}
\end{equation*}
So $\frac{n_1+n_3+n_5+n_7}{n_1n_2n_3}$ is a $3$-step relation number under the condition (\ref{3-step-1}). Hence we obtain the statement by choosing $n_1=r, n_2=s, n_3=t, n_4=1, n_5=v, n_6=1, n_7=w$.
\end{proof}
\begin{remark}
Proposition \ref{property-2step} and Proposition \ref{property-3step} imply that for arbitrary  nonzero integers $r,s,t,v, w$ satisfying either (i) 
$vw\,|\,rst$ and $rv+tv+tw=0$
or (ii) $rst=vw$ and $r^2s+w^2+rst=(r+t)(v+w)$, we have 
 $$-4<\frac{r+t+v+w}{rst}<4,$$
since any real number $q\in \mathbb R \setminus (-4,4)$ can not be a relation number (see \cite{Brenner} and \cite{Sanov}).
\end{remark}

\section{Proof of Theorem \ref{main-higher step}} \label{section:pf main-higher step}

Let $u$ be a nonzero root of $s^{\mathbf{n}}_{k+1}(t)$ for some integer sequence  ${\bf n}=(n_1,n_2,\cdots, n_{k})$, and let's assume that $s^{\mathbf{m}}_{l+1}(t)\neq 0$ for any integer sequence ${\bf m}$ if $l<k$. 
In this section, we show that such an integer sequence lies within a bounded region in $\mathbb{Z}^k$,  implying the existence of $N_i$, $i=1,2,\cdots,k$, such that $|n_i|\leq N_i$ for each $i$. 
 This proves Theorem \ref{main-higher step}  stated in the introduction.
For this purpose, we first consider the following well-known generalized continued fraction.


\begin{definition} \label{def-wcf}
    Let $\bar{\mathbb{C}}=\mathbb{C} \cup \{\infty\}$  and $\lambda$ be a nonzero complex number. A $\lambda$-$weighted$ $continued$ $fraction$ of length $k+1$ is a multivariate rational function from $\mathbb{Z}^{k+1}$ to $\bar{\mathbb{C}}$
    $$
        \wcf{1/\lambda}{{n_0; n_1,n_2,\cdots,n_k}} := n_0 + \cfrac{1/\lambda}{n_1 + \cfrac{1/\lambda}{n_2 + \cfrac{1/\lambda}{\cdots + \cfrac{1/\lambda}{n_k}}}} \in \bar{\mathbb{C}}.
    $$
\end{definition}
For an integer sequence, ${\bf n} = (n_1, n_2, \cdots, n_k) \in \mathbb{Z}^k$, we simply denote $\wcf{{1/\lambda}}{{0;n_1, \cdots, n_k}}$ by $\swcf{{\bf n}} = \swcf{{n_1, n_2, \cdots, n_k}}$. Let $\swcf{\ }:= 0$ for the empty sequence. For example, $\swcf{{\bf n}} = \swcf{n_1} = \frac{1 / \lambda}{n_1}$ when $k=1$.
\begin{proposition}\label{wcf Cheby}
Let ${\bf n} = (n_1, n_2, \cdots, n_k) \in \mathbb{Z}^k$. Then we have 
\begin{enumerate}
\item[\rm (i)]
    $\swcf{{n_1, n_{2}, \cdots, n_k}} = - \displaystyle\frac{s^{\sigma({\bf n})}_{k}(u)}{us^{{\bf n}}_{k+1}(u)}$, $\swcf{{n_k, n_{k-1}, \cdots, n_1}} = - \displaystyle\frac{s^{{\bf n}}_{k}(u)}{us^{{\bf n}}_{k+1}(u)}$
\item[\rm (ii)]
    $s^{\bf n}_{k+1}(u)=0$ if and only if $\swcf{{\bf n}} = \infty$.
    \end{enumerate}
\end{proposition}
\begin{proof}
We can prove the two identities in (i) by induction using the recursion
$$n_k u s_{k}^{{\bf n}}(u)-s_{k-1}^{\bf n}(u)=s_{k+1}^{\bf n}(u).$$
For example, for $k=1$,
$$\swcf{{\bf n}} = \swcf{n_1} = \frac{1 / \lambda}{n_1}=-\frac{s^{\bf n}_{1}(u)}{us^{{\bf n}}_{2}(u)}$$
and if we assume that the statement holds for all $k<l$, then the case $k=l$ is also proved by the induction hypothesis as follows.
\begin{equation*}
    \begin{split}
         \swcf{{n_l, n_{l-1}, \cdots, n_1}}& = \frac{1/\lambda}{n_l+ \swcf{{n_{l-1}, n_{l-2}, \cdots, n_1}}}\\
         &= \frac{-\frac{1}{u^2}}{n_l- \displaystyle\frac{s^{\bf n}_{l-1}(u)}{us^{{\bf n}}_{l}(u)}}= \frac{-\frac{1}{u}s^{\bf n}_{l}(u)}{n_lus^{\bf n}_{l}(u)-s^{\bf n}_{l-1}(u)}= -\frac{s^{\bf n}_{l}(u)}{us^{\bf n}_{l+1}(u)}.  \\
    \end{split}
\end{equation*}
This proves the second identity of (i), and the first identity can also be proved using the identity  $s_{k+1}^{\bf n}(u) = s_{k+1}^{\overleftarrow{\bf n}}(u)$  derived from (iii) of Lemma \ref{continuant-lemma1} and (i) of Lemma \ref{s-v-continuant}.
 
As we have seen earlier, $s_{k+1}^{\bf n}(u)$ and $s_k^{\bf n}(u)$ cannot be zero simultaneously due to the recursive definition of $s^{{\bf n}}$. Similarly, $s_{k+1}^{\bf n}(u)$ and $s_k^{\sigma({\bf n})}(u)$ cannot be zero simultaneously. Consequently, (ii) directly follows from (i).

     
\end{proof}

\begin{remark}
Kim-Koberda showed in \cite{KK} 
 that   a nonzero complex number $\lambda$ is a relation number if and only if there exists 
    ${\bf n} = (n_1, n_2, \cdots, n_k) \in (\mathbb{Z}\setminus \{0\})^k$
    such that 
    $\swcf{{n_k, n_{k-1}, \cdots, n_1}}=\infty$ for some $k\geq 2$ 
    or $\swcf{{n_k, n_{k-1}, \cdots, n_1}}=\swcf{{n_{k'}, n_{k'-1}, \cdots, n_1}}$ for some $k>k'\geq 2$. In fact,  in the latter case, it can be proven using Lemma \ref{ni-root-identity} and Proposition \ref{split s} that  
    $$\swcf{{n_1, n_2,\cdots, n_{k'-1}, n_{k'}-n_k, -n_{k-1}, \cdots, -n_1}}=\infty.$$
\end{remark}

Proposition \ref{split s} allows us to rewrite an equation $s^{\bf n}(u) = 0$ in terms of $\lambda$-weighted continued fractions as follows, which is needed to prove Theorem \ref{main-higher step}.


\begin{lemma}\label{root of Cheby and nseq}
For an integer sequence ${\bf n} = (n_1, n_2, \cdots, n_k)$ and for each $i=1, 2, \cdots, k$, $s^{\bf n}_{k+1}(u) = 0$ if and only if 
\begin{align*}
     n_i = -(\swcf{{n_{i-1}, \cdots, n_1}} + \swcf{{n_{i+1}, \cdots, n_{k}}})\, \text{ or }\, \swcf{{n_{i-1}, \cdots, n_1}} = \swcf{{n_{i+1}, \cdots, n_{k}}} = \infty.
\end{align*} Similarly,
\begin{align*}
    \tilde{T}^{\bf n}_{k}(u) = 0 \Leftrightarrow & \; n_i = - \left( \swcf{n_{i+1} \cdots n_k, n_1, \cdots, n_{i-1}} + \swcf{n_{i-1}, \cdots, n_1, n_k, \cdots, n_{i+1}} \right) \\
    & \text{ or }\, \swcf{n_{i+1} \cdots n_k, n_1, \cdots, n_{i-1}} = \swcf{n_{i-1}, \cdots, n_1, n_k, \cdots, n_{i+1}} = \infty
\end{align*}
and
\begin{align*}
    \hat{T}^{\bf n}_{k}(u) = 0 \Leftrightarrow \; & \,\, n_i = \swcf{n_{i+1} \cdots n_k, -n_1, \cdots, -n_{i-1}} - \swcf{-n_{i-1}, \cdots, -n_1, n_k, \cdots, n_{i+1}} \\
    &\text{or }\, \swcf{n_{i+1} \cdots n_k, -n_1, \cdots, -n_{i-1}} = \swcf{-n_{i-1}, \cdots, -n_1, n_k, \cdots, n_{i+1}} = \infty.
\end{align*}

\end{lemma}

\begin{proof}
    For each $i$, we have
    \begin{align}
        \label{s expand at i 3terms} s^{{\bf n}}(u) &= n_i u s^{{\bf n}_L}(u) s^{{\bf n}_R}(u) - s^{{\bf n}_{L'}}(u) s^{{\bf n}_R}(u) - s^{{\bf n}_L}(u) s^{{\bf n}_{R'}}(u) \\ 
        \label{s expand at i 2terms} &= s^{{\bf n}_L}(u) s^{(n_i, {\bf n}_R)}(u) - s^{{\bf n}_{L'}}(u) s^{{\bf n}_R}(u)        
    \end{align}
    by (i) of Proposition \ref{split s} or by the Euler identity (Proposition \ref{Euler Identity}). If none of $s^{{\bf n}_L}(u)$ and $s^{{\bf n}_R}(u)$ is zero, then $s^{{\bf n}}(u) = 0$ if and only if \begin{align*}
        n_i &= \frac{s^{{\bf n}_{L'}}(u)}{u s^{{\bf n}_L}(u)} + \frac{s^{{\bf n}_{R'}}(u)}{u s^{{\bf n}_R}(u)} \\
        &= -(\swcf{{n_{i-1}, \cdots, n_1}} + \swcf{{n_{i+1}, \cdots, n_{k}}})
    \end{align*}
    by \eqref{s expand at i 3terms}. For considering a case $s^{{\bf n}_L}(u)=0$ or $s^{{\bf n}_R}(u)=0$, note that \eqref{s expand at i 2terms} implies that if two of $s^{{\bf n}}(u)$, $s^{{\bf n}_L}(u)$ and $s^{{\bf n}_R}(u)$ are zero then the other must be zero. For example, if $s^{\bf n}(u)$ and $s^{{\bf n}_L}(u)$ are zero, then $s^{{\bf n}_{L'}}(u)$ cannot be zero by the definition of $s^{{\bf n}_L}(u)$; so it implies that $s^{{\bf n}_R}(u)=0$. This completes the proof for the statement for $s^{{\bf n}}(u)$.
    One can show the statements for $\tilde{T}^{{\bf n}}_k$ and $\hat{T}^{{\bf n}}_k$ by applying similar, but much simpler, arguments to \ref{split s split ttilde} and \ref{split s split that} of Proposition \ref{split s}. 
\end{proof}

For the proof of the Theorem \ref{main-higher step},
we need to extend the domain of $\lambda$-weighted continued fractions to include infinity as below.

\def\ms{\mathbf{m}}

\begin{lemma}\label{continued fraction infinity extension}
    Consider $\bar{\bz} = \bz \cup \{\infty\}$ as a subset of $\bar{\bc}$.
    Suppose $s_{l+1}^{\ms}(u) \neq 0$
    for any sequence $\ms \in (\bz \setminus \{0\})^l$ of length $l < k$.
    Then, ${\ns} \mapsto \swcf{\ns}$ is contiuously well-defined on $(\bar{\bz} \setminus \{0\})^k$. Furthermore, we have $[\nsL, \infty, \nsR] = [\nsL]$.
\end{lemma}
\begin{proof}
  We will show the following by induction on $k$.
  \begin{itemize}
      \item $[\ms] \neq \infty$ for a sequence $\ms$ of length $l < k$.
      \item $\ns \mapsto [\ns]$ is well-defined on $(\bar{\bz} \setminus \{0\})^k$.
      \item $[\nsL, \infty, \nsR] = [\nsL]$.
  \end{itemize}
  This is clear for $k = 1$.
  For $k > 1$, we first show $[\ms] \neq \infty$ for a sequence $\ms$ of length $l < k$.
  If $\ms$ does not contain infinity, this can be derived from $s_{l+1}^{\ms}(u) \neq 0$.
  Otherwise, one has
  \[
    [{\bf m}] = [{\bf m}_L, \infty, {\bf m}_R] = [{\bf m}_L] \neq \infty
  \]
  by induction hypothesis.
  This shows the well-definedness of $[\ns]$ as well, considering
  \[
    [\ns] = [n_1, \nsR] = [n_1 + [\nsR]]
  \]
  and $\nsR$ has length $k-1 < k$.
  Finally, we can check by computation that
  \[
    [\nsL, \infty, \nsR] = [\nsL, \infty + [\nsR]] = [\nsL]
  \]
  holds.
\end{proof}

\begin{proof}[\bf{Proof of Theorem \ref{main-higher step}}]
    It suffices to only consider sequences of nonzero integers, since we can replace ${\bf n}$ into a shorter sequence by using Lemma \ref{ni-root-identity} if ${\bf n}$ contains a zero. Suppose that $\ns$ is such a sequence which satisfies $s_{k+1}^\ns(u) = 0$ with minimal length.
    By Lemma \ref{root of Cheby and nseq}, we obtain
    \begin{align} \label{range ni}
      \begin{split}
        \abs{n_i} & = \abs{\swcf{{n_{i-1}, \cdots, n_1}} + \swcf{{n_{i+1}, \cdots, n_{k}}}} \\
        & \leq \max_{m_1, \cdots, m_{i-1} \in \bar{\bz} \setminus \{0\}} \abs{\swcf{{m_{i-1}, \cdots, m_1}}} + \max_{m_{i+1}, \cdots, m_k \in \bar{\bz} \setminus \{0\}} \abs{\swcf{{m_{i+1}, \cdots, m_{k}}}} < \infty.
      \end{split}
    \end{align}
    Consequently, we have a bound $N_i$ of the sequence where $\abs{n_i} \leq N_i$. Therefore, this implies that there could be only finitely many such sequences for a given $\lambda$ and $k$.
\end{proof}

Remark that similar approach could be applied to $\tilde{T}$, $\hat{T}$ and $T$ to obtain a bound of ${\bf n}$ satisfying $\tilde{T}^\ns(u) = 0$, etc.
In the case of $T$, $T^{(n_1, \cdots, n_k)} = 2 \: s^{(n_1, \cdots, n_{k-1}, n_k / 2)}$ by Lemma \ref{T-s identity}, so the same argument as for $s^{\bf n}_{k+1}$ applies to it. Bounds of such ${\bf n}$ for $\tilde{T}, \hat{T}$ can be obtained as follows.

\begin{proposition}
    An integer sequence of length $k$, ${\bf n}$, satisfying $\tilde{T}^{\ns}(u) = 0$ (or $\hat{T}^\ns(u) = 0$), lies within the range $\displaystyle \abs{n_i} \leq \max_{\mathbf{m} \in (\bar{\bz} \setminus \{0\})^{k-1}} \abs{\swcf{\mathbf{m}} + \swcf{\overleftarrow{\mathbf{m}}}}$(or $\displaystyle \max_{\mathbf{m} \in (\bar{\bz} \setminus \{0\})^{k-1}} \abs{\swcf{\mathbf{m}} - \swcf{\overleftarrow{\mathbf{m}}}}$, respectively), provided that $s_{l+1}^\ms(u) \neq 0$ for any sequence $\ms$ of length $l < k$.
\end{proposition}
\begin{proof}
    As in the proof of Theorem \ref{main-higher step}, we have the $\lambda$-weighted continued fraction well-defined on $(\bar{\bz} \setminus \{0\})^k$,
    with corresponding maximum for each length $l \leq k$. Recall that in  Lemma \ref{root of Cheby and nseq},
    each such a sequence satisfies these conditions:
    \begin{align*}
      \tilde{T}^{\ns}(u) = 0 \, \Leftrightarrow &
      \quad n_i = - \swcf{\nsR, \nsL} - \swcf{\overleftarrow{\nsL}, \overleftarrow{\nsR}}, \\
      \hat{T}^\ns(u) = 0 \, \Leftrightarrow &
      \quad n_i = \swcf{\nsR, -\nsL} - \swcf{- \overleftarrow{\nsL}, \overleftarrow{\nsR}}.
    \end{align*} In both cases, we obtain
    \[
      \abs{n_i}
      \leq \max_{\nsL, \nsR} \abs{\mp \swcf{\nsR, \pm \nsL} - \swcf{\pm \overleftarrow{\nsL}, \overleftarrow{\nsR}}}
    \]
    from which the desired conclusion follows.
\end{proof}

\section{Algorithms for relation numbers} \label{section:algorithm}

In this section, we prove Theorem \ref{decision algorithm} by establishing an exhaustive algorithm to find a nonzero integer sequence ${\bf n}=(n_1, n_2, \cdots, n_k)$ of a given length $k$ satisfying $s^{\bf n}_{k+1}(\sqrt{-\lambda}) = 0$ for a given $\lambda \in \mathbb{C}$ assuming $k$ is the minimal such length. Using the fact that  $s^{\bf n}_{k+1}(u)=0$ if and only if $\swcf{{\bf n}} = \infty$ (Proposition \ref{wcf Cheby}), 
our exhaustive algorithm is designed to  find the shortest nonzero integer sequence ${\bf n}=(n_1, n_2, \cdots, n_k)$ satisfying $\swcf{{\bf n}}=\infty$ 
 by inductively computing the maximum absolute value $M_i$ of $\lambda$-weighted continued fractions for each increasing length $i$ and checking whether $M_i$ equals to $\infty$ or not.

We first introduce an exhaustive algorithm for real $\lambda \in \mathbb{R}$, and then we discuss an extension of the algorithm to the complex case, i.e., for  $\lambda \in \bc$, with a few modifications.
We also present a modified version, called the greedy algorithm, designed for practical purpose for finding an ${\bf n}$-sequence more easily, but not necessarily shortest. 

For this purpose, we introduce the following terminologies.
\begin{definition}
    For $\lambda=-u^2 \in \mathbb{C}$, the \textit{minimal $u$-degree} of $\lambda$, $\mindeg_u(\lambda)$, is the minimal length $k$ of a nonzero integer sequence ${\bf n}=(n_1, n_2, \cdots, n_k)$ such that $s^{{\bf n}}_{k+1}(u)=0$. If there is no such ${{\bf n}}$, then $\mindeg_u(\lambda) := \infty$.
\end{definition}



\begin{definition}
    For a fixed nonzero complex $\lambda$ and a positive integer $k$, the $k$-$th$ $maximum$ $absolute$ $value$ of $\lambda$-weighted continued fraction, valued in $\mathbb{R}_{>0} \cup \{\infty\}$, is defined as
    $$
        \wcfM{k}(\lambda) := \begin{cases}
            \infty   &\text{ if } s_{j+1}^\ns(u) = 0 \text{ for some } {\bf n} \text{ of length } j < k, \\
           \max{\{ \abs{\swcf{{\bf n}}} \,|\, {\bf n} \in (\bar{\mathbb{Z}}\setminus \{0\})^k\}} &\text{ otherwise,}
        \end{cases}
    $$
     and let $\wcfM{0}(\lambda) := 0$.
\end{definition}
 Note that the minimum $k$ such that $M_k(\lambda) = \infty$ is the same as $\mindeg_u(\lambda)$. We use the shorthand notation $M_k$ for $M_k(\lambda)$ whenever $\lambda$ is well-understood.

Let $k_0 = \mindeg_u(\lambda)$. Consider a sequence \begin{align} \label{Mk sequence}
    (0=M_0, M_1, M_2, \cdots, M_{k_0-1}, M_{k_0}=\infty)
\end{align}
which might be an infinite sequence when $k_0=\infty$. By Lemma \ref{continued fraction infinity extension}, $\swcf{{n_1, n_2, \cdots, n_k}} = \swcf{{n_1, n_2, \cdots, n_k, \infty}} \leq M_{k+1}$ for $k < k_0$, and that the above sequence is an increasing sequence,
i.e., if $k_0 < \infty$
\[
    0 = M_0 < M_1 \leq M_2 \leq \cdots \leq M_{k_0 - 1} < M_{k_0} = \infty.
\]

Remark that $k_0 < \infty$ if and only if $\lambda$ is a relation number. For $\lambda \in \mathbb{R}$, it suffices to consider only positive $\lambda$ since $M_k(\lambda) = M_k(-\lambda)$. Note that \begin{align} \label{wcf sign change}
\wcf{-1/\lambda}{{n_1, n_2, n_3, \cdots, n_k}} = \wcf{1/\lambda}{{-n_1, n_2, -n_3, \cdots, (-1)^k n_k}}    
\end{align} from the definition. The following proposition offers a more simplified understanding of the sequence \eqref{Mk sequence} when $\lambda$ is a real number.

\begin{proposition} \label{Mk only nonzero integer}
    Let $k_0 \in \mathbb{N}\cup\{\infty\}$ be the minimal $u$-degree of a nonzero real $\lambda$. Then, for all positive integer $k \leq k_0$,
    $$
    M_k = \max\{\abs{\swcf{{n_1, n_2, \cdots, n_k}}} \,|\, n_1, n_2, \cdots, n_k \in \mathbb{Z}\setminus \{0\}\},$$
    i.e., $M_k$ is attained with $n_i \neq \infty$ for all $i=1,2,\cdots,k$, and furthermore $M_{k-1} < M_k$.    
    
\end{proposition}
\begin{proof}
Fix a nonzero real $\lambda$ and $u \in \sqrt{-1}\; \mathbb{R}$ such that $\lambda = -u^2$.
We may consider $\swcf{{\bf n}}$ as a rational function from $\br^k$ to $\bar{\br} = \br \cup \{\infty\}$;
$\swcf{{\bf n}}$ can be written as $t(n_1, t(n_2, \cdots t(n_k, 0) \cdots)$ where $t(z, w) = \frac{1}{\lambda} \cdot \frac{1}{z+w}$ and $\pdv{}{z} t(z,w) = \pdv{}{w} t(z,w) = -\lambda (t(z, w))^2$. One can directly calculate 
\begin{align*}
    \pdv{}{n_i} \swcf{\bf n} 
    &= \prod_{j=1}^i (-\lambda \swcf{{n_j, n_{j+1}, \cdots, n_k}}^2) \neq 0.
\end{align*}
The fact that the above expression is not zero arises from each term $\swcf{n_j, n_{j+1}, \cdots, n_k}$ being nonzero, which is ensured by Proposition \ref{wcf Cheby} and the assumption of minimal $u$-degree, ensuring that $s_{k-j+1}^{(n_{j+1}, \cdots, n_k)}(u) \neq 0$.
 Moreover, by Proposition \ref{split s},
\begin{align} \label{ni rational}
\swcf{{\bf n}} &= -\frac{1}{u} \cdot \frac{n_iu s^{{\sigma({\bf n}_L)}}_{i-1}(u) s^{{\bf n}_R}_{k-i+1}(u) - (s^{{\sigma({\bf n}_{L'})}}_{i-2}(u) s^{{\bf n}_R}_{k-i+1}(u) + s^{{\sigma({\bf n}_L)}}_{i-1}(u) s^{{\bf n}_{R'}}_{k-i}(u)) }{n_iu s^{{\bf n}_L}_{i}(u) s^{{\bf n}_R}_{k-i+1}(u) - (s^{{\bf n}_{L'}}_{i-1}(u) s^{{\bf n}_R}_{k-i+1}(u) + s^{{\bf n}_L}_{i}(u) s^{{\bf n}_{R'}}_{k-i}(u)) } \quad
\\ \nonumber &
\text{  (note that $s^{\bf n}_{-1}(u) = -1$)}
\\ 
\nonumber &= \begin{cases}
\displaystyle \frac{1 / \lambda}{n_1 + [n_2, \cdots, n_k]}, & i = 1 \\
  \displaystyle \swcf{{n_{1}, \cdots, n_{i-1}}} \cdot \frac{n_i + \swcf{{n_{i-1}, \cdots, n_2}} + \swcf{{n_{i+1}, \cdots, n_k}}}{n_i + \swcf{{n_{i-1}, \cdots, n_1}} + \swcf{{n_{i+1}, \cdots, n_k}}}, & i \geq 2.
\end{cases}
\end{align} In other words, $\swcf{{\bf n}}$ is a rational function of $n_i=x$ with numerator and denominator both polynomials of degree $\leq 1$, which has a singular point
\begin{align} \label{singular point}
x_\infty &= -(\swcf{{n_{i-1}, \cdots, n_1}} + \swcf{{n_{i+1}, \cdots, n_{k}}})
\end{align}
In particular, one can always choose one of nearest integers to the singular point (which might be the point itself), say $m_i$, to maximize the absolute value with $\abs{\swcf{{n_1, \cdots, m_i, \cdots, n_k}}} \geq \abs{\swcf{{n_1, \cdots, n_i, \cdots, n_k}}}$ for all $n_i \in \bz \setminus \{0\}$. See Figure \ref{fig:[n]graph} below.

\begin{figure}[h]
    \centering
        \begin{tikzpicture}
            \begin{axis}[
                axis x line=middle,
                axis y line=middle,
                xmin=-4, xmax=2,                
                ymin=-8, ymax=12.5,
                ticks=none,
                x label style={at={(ticklabel* cs:1.)},anchor=west},
                xlabel = $x$,
                y label style={at={(ticklabel* cs:1.)},anchor=south},
                ylabel = $y$,
                width=.6\textwidth,
                clip mode=individual,
            ]
            \addplot[
                samples=1001,
                very thick,
                unbounded coords=jump,                
            ]{(2*x+1)/(x+1)};  
            \addplot[mark=none, dashed] coordinates {(-6, 2) (4, 2)};
            \addplot[mark=none, dashed] coordinates {(-1,-8) (-1,12.5)};
            \addplot[mark=none, dotted] coordinates {(-1.2,0) (-1.2,6.9)};
            \node[label={[label distance=-3pt]225:{$m_i$}},circle,fill,inner sep=1.2pt] at (axis cs:-1.2,0) {};
            \node[label={[shift={(24pt,15pt)}]225:{$x_\infty$}},circle,fill,inner sep=1.2pt] at (axis cs:-1,0) {};
            \node[label={[shift={(-36pt,12pt)}]0:{$y=\swcf{{n_1, \cdots, n_{i-1}}}$}}] at (axis cs:2,2) {};
            \node[label={[label distance=-8pt]-45:$O$}] at (axis cs:0,0) {};
            \node at (axis cs:-2,6) {$y=\swcf{{\bf n}}$};            
            \end{axis}
        \end{tikzpicture}
    \caption{Graph of $y=\swcf{{\bf n}}$ as a function of $x=n_i$ along with its asymptotic line $y = \swcf{{n_1, \cdots, n_{i-1}}}$, when $i$ is even.}
    \label{fig:[n]graph}
\end{figure}
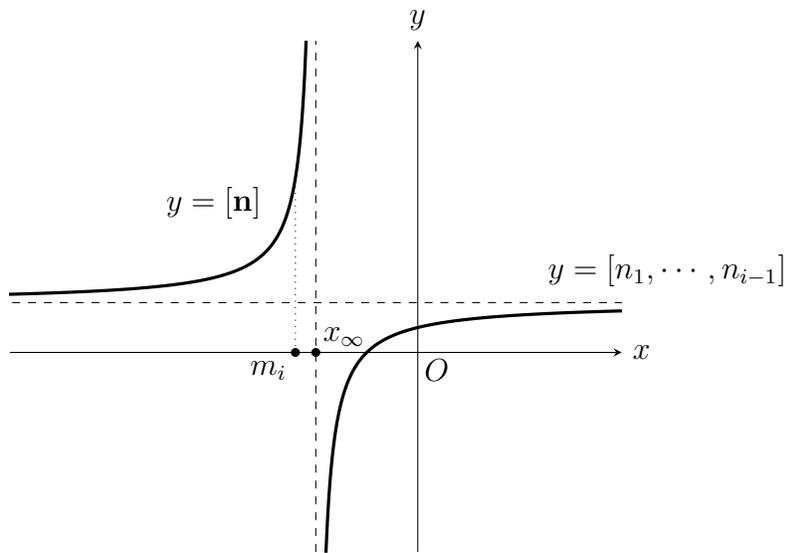

Taking $i = k+1$, one has $m_{k+1} \in \bz \setminus \{0\}$ 
satisfying
\begin{equation}\label{has bigger fraction}
    \abs{\swcf{n_1, \cdots, n_k}} < \abs{\swcf{n_1, \cdots, n_k, m_{k+1}}}
\end{equation} considering the behavior when $\abs{n_{k+1}} \to \infty$, with other $n_j$'s being finite.
Suppose that $M_k$ is attained at ${\bf n}$ and let say $n_i$ is the leftmost argument which is $\infty$. By Proposition \ref{continued fraction infinity extension} and the above argument, we have $m_{i+1} \in \bz \setminus \{0\}$ such that
\[
M_k =  \swcf{n_1, \cdots, n_{i-1}, \infty, n_{i+1}, \cdots, n_k} = \swcf{{n_1, \cdots, n_{i-1}}} < \swcf{{n_1, \cdots, n_{i-1}, m_i}} \leq M_i,
\]
leading to a contradiction. Therefore, $M_k$ is always attained at a finite integer sequence of length $k$.
Lastly, \eqref{has bigger fraction} shows that $M_k$ is a strictly increasing sequence.
\end{proof}

Note that a singular point of a rational function \eqref{ni rational} lies in the range \eqref{range ni} and the nearest integer points to this singular point are candidates for maximizing $\abs{\swcf{{\bf n}}}$.
\begin{proposition} \label{max wcf}
    Let $\lambda$ be a nonzero real number and $k_0$ be its minimal $u$-degree. For $k \leq k_0$, each $\wcfM{k}$ is attained by a nonzero integer sequence ${\bf n}=(n_1, n_2, \cdots, n_k)$ satisfying
    \[
    \abs{n_i + \swcf{{n_{i-1}, \cdots, n_1}} + \swcf{{n_{i+1}, \cdots, n_k}}} < 2.
    \]
    In particular,
    \begin{align*}
        \abs{n_i + \swcf{{n_{i-1}, \cdots, n_1}}} & < M_{k-i} + 2
    \end{align*}
    for all $i=1, 2, \cdots, k$.
\end{proposition}

\begin{proof}

    Let us apply induction on $k$.
    Note that
    $$
        \abs{\swcf{n_1}} = \abs{\frac{1}{\lambda n_1}} \leq \frac{1}{\abs{\lambda}} = \wcfM{1}
    $$ for all $n_1 \in \mathbb{Z} \setminus \{0\}$. Hence, $\wcfM{1}$ is obtained when $|n_1| = 1 < \wcfM{0} + 2$.

    It suffices to prove for $k > 1$.
    We will consider $[\ns]$ to be a function of $n_i$ only, and fix all the other $n_j$'s.
    To emphasize this, denote $f(x) := \swcf{{n_1, \cdots, n_{i-1}, x, n_{i+1}, \cdots, n_k}}$.
   Note that $x_{\infty}$ is the singular point of $f(x)$. 
    Then, take $l_i$, $r_i$ to be the floor and ceiling of $x_\infty$, respectively.
    If the singular point $x_\infty$ is a nonzero integer, one has $M_k = \infty$ obtained from the nonzero integer sequence $(n_1, \cdots, n_k)$ with $n_i = x_\infty$.
    Otherwise, $l_i < x_\infty < r_i$ and so
    $$
        \begin{cases}
            f(r_i) \leq f(m) \leq f(l_i) &\text{if } i \text{ is even } \\
            f(l_i) \leq f(m) \leq f(r_i) &\text{if } i \text{ is odd }
        \end{cases} \Rightarrow |f(m)| \leq \max{\{ \abs{f(r_i)}, \abs{f(l_i)} \}}
    $$ for all $m \in \mathbb{Z} \setminus \{0\}$. If one of $l_i$ or $r_i$ equals to zero, then one may take $l_i=-1$ and $r_i=1$ instead without violating any of the inequalities above. As a consequence, if $\wcfM{k}$ is obtained by ${\bf n} = (n_1, n_2, \cdots, n_{k}$) then nonzero $n_i$ must be either $l_i$ or $r_i$ for each $i = 1, 2, \cdots k$.
\end{proof}

Following the proof of Proposition \ref{max wcf}, we can calculate a sequence $(0=M_0, M_1, M_2, \cdots)$ in a stepwise manner. Start with $M_1 = \frac{1}{\abs{\lambda}}$, then determine $M_2$ by computing the maximum of $\abs{\swcf{{n_1, n_2}}}$ among nonzero integer sequences ${\bf n}$ within the range $\abs{n_1} < M_1 + 2$ and $\abs{n_2 + [n_1]} < M_0 + 2 = 2$. If $M_2 = \infty$ is obtained by one of the sequences, then the sequence should satisfy $s^{\bf n}_3(u) = 0$. Otherwise, find $M_3$ by calculating the maximum of $\abs{\swcf{{n_1, n_2, n_3}}}$ within a range $\abs{n_1} <  M_2 + 2$, $\abs{n_2 + [n_1]} < M_1 + 2$, and $\abs{n_3 + [n_2, n_1]} < M_0 + 2$, and repeat this procedure. Fortunately, one might skip the initial steps by the following proposition.

\begin{proposition} \label{max wcf alternating test}
    Let $\lambda$ be a fixed nonzero real number. If $\wcfM{k} \leq 1$ for some positive integer $k$, then 
    $$
    \wcfM{k+1} = \begin{cases}
        \swcf{{1, -1, 1, \cdots, (-1)^k}} & \text{ if } \lambda > 0, \\
        \swcf{{-1, -1, -1, \cdots, -1}} & \text{ if } \lambda < 0.
    \end{cases}
    $$
\end{proposition}

\begin{proof}
Let $\lambda = -u^2 \in \mathbb{R}_{>0}$ and let us prove the statement by induction. First, note that $M_1 = \swcf{{1}}$ as checked in the proof of Proposition \ref{max wcf}. Assume the statement holds for $k-1$ and $M_{k} \leq 1$. Then, by Proposition \ref{Mk only nonzero integer}, $M_{k-1} < M_{k} \leq 1$ so that $M_{k} = \swcf{{1, -1, 1, \cdots, (-1)^{k-1}}}$ by the induction hypothesis. Since $\abs{\swcf{n_2, \cdots, n_{k+1}}} \leq 1$ for all nonzero integer sequence ${\bf n}$, 
$$
\swcf{{n_1, n_2, \cdots, n_{k+1}}} = \frac{1}{\lambda} \cdot \frac{1}{n_1 + \swcf{n_2, \cdots, n_{k+1}}}
$$
is maximized when $n_1 = 1$ and $\swcf{n_2, \cdots, n_{k+1}}$ is closest to $-1$ which is the case that $\swcf{n_2, \cdots, n_{k+1}} = -M_k = \swcf{-1, 1, \cdots, (-1)^k}$ by the induction hypothesis (note that $-\swcf{{\bf n}} = \swcf{-{\bf n}}$). Therefore, $M_{k+1} = \swcf{{1, -1, \cdots, (-1)^k}}$. The statement for $\lambda < 0$ follows directly from \eqref{wcf sign change}.
\end{proof}

Now we describe a structure for an algorithm which checks whether a given positive real $\lambda$ is a relation number or not by combining Proposition \ref{max wcf} and \ref{max wcf alternating test}:
First, calculate $\swcf{{1, -1, \cdots, (-1)^k}}$ until it exceeds $1$. Suppose that $M_{k'-1} \leq 1 < M_{k'} = \swcf{{1, -1, \cdots, (-1)^{k'-1}}}$ for some $k' > 0$. Let $k=k'+1$ and calculate $\abs{\swcf{{n_1, n_2, \cdots, n_k}}}$ for all ${{\bf n}} = (n_1, n_2, \cdots, n_k)$ satisfying
\begin{align*}
    \abs{n_1} & < M_{k-1} + 2\\
    \abs{n_2 + \swcf{{n_1}}} & < M_{k-2} + 2\\
    \abs{n_3 + \swcf{{n_2, n_1}}} & < M_{k-3} + 2\\
    \ \vdots & \\
    \abs{n_k + \swcf{{n_{k-1}, n_{k-2}, \cdots, n_1}}} & < M_{0} + 2 = 2.
\end{align*}
If there exists a sequence with $\swcf{{\bf n}} = \infty$ in the above range, then we have $M_k = \infty$ and the algorithm terminates. Otherwise, let $M_{k}$ be the maximum value among calculations, and repeat the process after increasing $k$ by $1$. See Algorithm \ref{max wcf alg} for the pseudocode. 

\begin{algorithm}[!ht]
  \caption{Exhaustive algorithm (for $\lambda > 0$)} \label{max wcf alg}
  \begin{algorithmic}[1]
    \State{$k \gets 0$, $M_0 \gets 0$}
    \While{$M_k \leq 1$}
        \State{${\bf n} \gets (+1, -1, \cdots, (-1)^k)$}
        \State{$M_{k+1} \gets \swcf{{\bf n}}$, $k \gets k+1$}
    \EndWhile   

    \If{$M_k = \infty$}
        \State {\Return ${\bf n}$}
    \Else{
        \While{$M_k < \infty$}
            \For{$n_1 \in (\mathbb{Z} \setminus \{0\}) \cap \left(-(M_{k-1}+2), (M_{k-1}+2)\right)$}
                \For{$n_2 \in (\mathbb{Z} \setminus \{0\}) \cap \left(-\swcf{{n_1}}-(M_{k-2}+2)\right.$, $\left. -\swcf{{n_1}}+(M_{k-2}+2)\right)$}
                    \State{$\vdots$}
                    \For{$n_{k-1}\in(\mathbb{Z} \setminus \{0\}) \cap \left(-\swcf{{n_{k-2},\cdots,n_1}}-(M_1+2)\right.$, $\left.-\swcf{{n_{k-2}, \cdots, n_1}}+(M_{1}+2)\right)$}
                        \State{$m \gets -\swcf{{n_{k-1}, \cdots, n_1}}$}
                        \If{$m \in \mathbb{Z}\setminus \{0\}$}
                            \State{$n_k \gets m$}
                            \State{\Return ${\bf n} = (n_1, n_2, \cdots, n_k)$} 
                        \Else{
                            \State{$l \gets \lfloor m \rfloor$ and $r \gets \lceil m \rceil$}
                            \If{$l = 0$ or $r = 0$}
                                \State $l \gets -1$, $r \gets 1$
                            \EndIf
                            \State{$M_k \gets \max{\{M_k, \abs{\swcf{{n_1, \cdots, n_{k-1}, l}}}, \abs{\swcf{{n_1, \cdots, n_{k-1}, r}}\}}}$}
                        }
                        \EndIf
                    \EndFor
                \EndFor
            \EndFor
        \EndWhile
    }
    \EndIf
  \end{algorithmic}
\end{algorithm}

Unfortunately, reliance of Algorithm \ref{max wcf alg} on nested $k$-loops becomes a serious drawback when the algorithm slows down considerably for large $k$. In particular, we experimentally observed an exponential growth tendency in $M_k$ as $k$ increases. For example, a sequence of maximum absolute values for $\lambda = \frac{27}{8}$ is as follows.
\begin{table}[H]
    \centering
\begin{tabular}{c|c|c|c|c|c|c|c|c|c|c|c}\hline
  $k$ & $1$ & $2$ & $3$ & $4$ & $5$ & $6$ & $7$ & $8$ & $9$ & $10$ & $\cdots$ \\ \hline
  \multirow{3}{2em}{\centering $M_k$} & & & & & & & & & & & \\
  & $\displaystyle \frac{8}{27}$ & $\displaystyle \frac{8}{19}$ & $\displaystyle \frac{152}{297}$ & $\displaystyle \frac{88}{145}$ & $\displaystyle \frac{1160}{1539}$ & $\displaystyle \frac{456}{379}$ & $\displaystyle \frac{3032}{2079}$ & $\displaystyle \frac{6680}{953}$ & $\displaystyle \frac{140680}{243}$ & $\displaystyle \frac{563336}{19}$ & $\cdots$ \\
  & & & & & & & & & & & \\ \hline
\end{tabular}
    \vspace*{3mm}
    \caption{A maximum absolute values of $27/8$-weighted continued fractions}
    \vspace*{-3mm}
    \label{tab:maxabs_wcf_27/8}
\end{table}
In this context, $\mindeg_u(\lambda)$ reflects difficulty of deciding whether a given $\lambda$ is a relation number or not. The following statement, which is a generalization of the well-known fact that $\lambda \geq 4$ cannot be a relation number, illustrates that the difficulty increases as $\lambda$ gets closer to $4$.

\begin{proposition} \label{lower bound min deg}
If $\lambda$ is a real number such that $\abs{\lambda} \geq 4\cos^2(\pi / (n+1))$ for some $n > 1$, then the minimal $u$-degree of $\lambda$ is greater than or equal to $n$. The equality holds only if $\lambda = \pm 4\cos^2(\pi/(n+1))$. In particular, $4 \cos^2(\pi/(2l + 2))$ is the largest $l$-step relation number among real numbers.

\end{proposition}
\begin{proof}
    It is enough to consider for $\lambda > 0$. By Lemma \ref{T_alt}, one can calculate that, for all $k>0$, \begin{align*}
    \swcf{{1, -1, 1, \cdots, (-1)^{k-1}}} &= -\frac{s^{{\bf n}}_{k}(u)}{us^{{\bf n}}_{k+1}(u)} = \lambda^{\frac{(-1)^{k}-1}{2}}\frac{q_k(\lambda)}{q_{k+1}(\lambda)}
\end{align*} 
where ${\bf n} = (1, -1, 1, -1, \cdots)$ and $q_k(\alpha)$ is a monic polynomial whose roots are $4\cos^2(\frac{1}{k} \pi)$, $4 \cos^2(\frac{2}{k} \pi)$, $\cdots$, $4 \cos^2(\frac{\lfloor k/2 \rfloor}{k} \pi)$ except $4 \cos^2(\pi/2)(=0)$. With the assumption of the statement, we have $F_k:=\swcf{{1, -1, 1, \cdots, (-1)^{k-1}}} > 0$ and $ 0 < F_{k-1} < 1$ for $k=2, 3, \cdots, n$ since $F_k = \swcf{{1 - F_{k-1}}} > 0$. Finally, Proposition \ref{max wcf alternating test} implies that $M_n (= \swcf{{1, -1, 1, \cdots, (-1)^{n-1}}})$ is $\infty$ only if $\lambda = 4\cos^2(\pi/(n+1))$.
\end{proof}

Thereby, a rational number $(4q-1)/q$ poses a significant challenge when deciding whether it is a relation number or not. In our experiments, Algorithm \ref{max wcf alg} terminates within a few minutes for $q=1, 2, 3$.
However, for $q > 3$ the Algorithm failed to produce results for hours.
We dealt with such cases by applying heuristics to obtain a `greedy algorithm'.

\begin{algorithm}[!ht]
  \caption{Greedy algorithm} \label{wcf greedy alg}
  \begin{algorithmic}[1]
    \Function{ContinuantZero}{$\lambda$, $N$, $W$, $k$}
        \For{$n_1 \in (\mathbb{Z} \setminus \{0\}) \cap  [-W, W]$}
            \For{$n_2 \in (\mathbb{Z} \setminus \{0\}) \cap [-\wcf{1/\lambda}{{n_1}}-W, -\wcf{1/\lambda}{{n_1}}+W]$}
                \For{$n_3 \in (\mathbb{Z} \setminus \{0\}) \cap [-\wcf{1/\lambda}{{n_2,n_1}}-W, -\wcf{1/\lambda}{{n_2,n_1}}+W]$}
                    \State{$\vdots$}
                    \For{$n_N \in (\mathbb{Z} \setminus \{0\}) \cap [-\wcf{1/\lambda}{{n_{N-1},\cdots,n_1}}-W, -\wcf{1/\lambda}{{n_{N-1}, \cdots, n_1}}+W]$}
                    \For{$i \in [N+1, k]$}
                        \State{$m \gets -\wcf{1/\lambda}{{n_{i-1}, \cdots, n_1}}$}
                        \State{$n_i \gets$ \Call{Round}{$m$}}
                        \If{$n_i = 0$}
                            \State $n_i \gets 1$
                        \EndIf
                        \If{$m \in \mathbb{Z}\setminus \{0\}$}
                            \State{\Return $\ns = (n_1, n_2, \cdots, n_k)$} \label{max wcf alg found infty}
                        \EndIf
                    \EndFor
                    \EndFor
                \EndFor
            \EndFor
        \EndFor
        \State{\Return fail}
    \EndFunction
  \end{algorithmic}
\end{algorithm}

The greedy algorithm selects small search range with the intuition earned in the exhaustive algorithm. Notice that in Algorithm \ref{max wcf alg}, the search range of $n_i \approx \swcf{{n_1, \cdots, n_{i-1}}} + \swcf{{n_k, \cdots, n_{i+1}}}$ is based around the previous continued fraction $\swcf{{n_1, \cdots, n_{i-1}}}$. Moreover, the right-hand side continued fraction $\swcf{{n_k, \cdots, n_{i+1}}}$ is found to be often smaller than $1$. This phenomenon suggests a heuristic of taking narrow search range around this point.
To be precise, we consider the following range as in Algorithm \ref{wcf greedy alg}:
\[
\begin{cases}
    \abs{n_i - \swcf{{n_1, \cdots, n_{i-1}}}} \leq W_\lambda, & 1 \leq i \leq N_\lambda \\
    n_i = \operatorname{round}(\swcf{{n_1, \cdots, n_{i-1}}}), & N_\lambda < i \leq k. \\
\end{cases}
\]
The values of $W_\lambda$ and $N_\lambda$ are taken through experiments. In the case of rational $\lambda = p/q$, we found that $W_{p/q} = q$ and $N_{p/q}=5 \text{ or } 6$ works well.
With such an appropriate $W_\lambda$ and $N_\lambda$, this algorithm yields ${\bf n}$-sequences for most of small rational $\lambda$ values. 

Our results are listed in the form of an additional file for $q < 23$, which gives the desired sequences for all such $p/q$ except for
\[
\frac{51}{13}, 
\frac{67}{17},
\frac{74}{19}, \frac{75}{19},
\frac{83}{21}, 
\frac{87}{22}
\]
and the result for $q=13$ is given in Table \ref{tab:result for q=13} as a sample. For each $\lambda = \frac{p}{13}$, we keep a record of ${\bf n}$ such that $s^{{\bf n}}_{k+1}(u) = 0$, and we designate $k$, the length of ${\bf n}$, as the minimal $u$-degree of $\lambda$ if ${\bf n}$ is found by Algorithm \ref{max wcf alg}. Otherwise, $k$ serves as an upper bound for the minimal $u$-degree. The lower bound is given by Proposition \ref{lower bound min deg}.

\begin{table}[!htp]
    \centering
    \begin{adjustbox}{height=.5\textheight, center}
        \begin{tabular}{c|c|l}
             $\lambda\,(=-u^2)$ & $\,\,d:=\mindeg_u(\lambda)\,\,$
              & ${\bf n} = (n_1, n_2, \cdots, n_k)$ such that $s^{{\bf n}}_{k+1}(u) = 0$ \\
             \hline
             1/13 & 2 & $(-13, 1)$\\
             2/13 & 3 & $(13, -7, 1)$\\
             3/13 & 3 & $(13, 4, -1)$\\
             4/13 & 3 & $(13, 3, -1)$\\
             5/13 & 5 & $(-13, 8, -9, -1, 2)$\\
             6/13 & 3 & $(13, 2, -1)$\\
             7/13 & 3 & $(13, -2, 1)$\\
             8/13 & 5 & $(13, -5, 9, -1, 2)$\\
             9/13 & 5 & $(13, -3, 1, -3, 13)$\\
             10/13 & 5 & $(-13, 4, 4, 1, -1)$\\
             11/13 & 5 & $(1, -1 ,-13, -1, 1)$\\
             12/13 & 3 & $(13, 1, -1)$\\
             14/13 & 3 & $(13, -1, 1)$\\
             15/13 & 5 & $(1, -1, 13, -1, 1)$\\
             16/13 & 5 & $(13, 4, -1, 1, 52)$\\
             17/13 & 5 & $(13, 29, 1, -1, 3)$\\
             18/13 & 7 & $(1, -1 ,3, -1 -1, 7, -39)$\\
             19/13 & 5 & $(13, 2, -1, 1, 26)$\\
             20/13 & 5 & $(13, -2, 1, -1, 26)$\\
             21/13 & 5 & $(-39, 6, -1, 1, -2)$\\
             22/13 & 7 & $(2, 26, 1, -1, 2, -1, 39)$\\
             23/13 & 7 & $(1, 65, 1, -1, 1, 2, -26)$\\
             24/13 & 5 & $(-26, 3, 1, -1, 1)$\\
             25/13 & 5 & $(13, 1, -1, 1, 13)$\\
             27/13 & 5 & $(13, -1, 1, -1, 13)$\\
             28/13 & 7 & $(1, -13, 1, 3, -1, 1, -1)$\\
             29/13 & 9 & $(1, -2, 13, -3 ,1, -1, 1, -9, 13)$\\
             30/13 & 7 & $(13, 3, -2, 1, -1, 1, -130)$\\
             31/13 & 9 & $(3, -13, 3, -1, 2, -1, 1, -1, 611)$\\
             32/13 & 9 & $(2, -299, 1, -1, 1, -1, -1, -1, -195)$\\
             33/13 & 9 & $(13, -2, 9, 1, -1, 1, -1, -6, -52)$\\
             34/13 & 5 & $(104, 1, -1, 1, -1)$\\
             35/13 & 7 & $(52, -4, 1, -1, 1, -1, -299)$\\
             36/13 & 7 & $(104, 19, 1, -1, 1, -1, 2)$\\
             37/13 & 9 & $(13, 7, -2, 1, -1, 1, -1, 4, 65)$\\
    
             38/13 & $ 12 \leq d \leq 151 $ & $(1,1,1, \cdots, 43,1,13)$ \\
             40/13 & $ 11 $ & $(2,-6,-11, \cdots, -3,-1040,-1346)$ \\
             41/13 & $ 13 \leq d \leq 125 $ & $(1,1,1, \cdots, 11,-6,-13)$ \\
             
             42/13 & 9 & $(13, 1, -1, 1, -1, 1, -1, -13, -94)$\\

             43/13 & $ 13 \leq d \leq 127 $ & $(1, 1, 1, \cdots, 1, 1, 39)$\\
             44/13 & $ 14 \leq d \leq 165 $ & $(1, 1, 1, \cdots, 1,-8,13)$ \\
             45/13 & $ 14 \leq d \leq 175 $ & $(1, 1, 1, \cdots, -1,78,-3)$ \\
             46/13 & $ 15 \leq d \leq 125 $ & $(1, 1, 1, \cdots, 17,-1,26)$\\
             47/13 & $ 15 \leq d \leq 107 $ & $(1, 1, 1, \cdots, -25,-13,2)$\\
             48/13 & $ 17 \leq d \leq 117 $ & $(1, 1, 1, \cdots, -7,13,-2)$ \\
             49/13 & $ 18 \leq d \leq 95 $ & $(1, 1, -1, \cdots, 320,13,-7)$ \\
             50/13 & $ 21 \leq d \leq 147 $ & $(1, 1, 7, \cdots, -4,-5,-52)$ \\
             51/13 & $ 29 \leq d $ & Unknown \\
             \hline
             
        \end{tabular}
    \end{adjustbox}
    \caption{A combined result of the exhaustive and greedy algorithm for rational $\lambda$ whose denominator is $13$}
    \label{tab:result for q=13}
\end{table}

\subsection*{Exhaustive algorithm for $\lambda \in \bc$}
\phantom{ }

The process of finding maximum of the continued fraction can be extended to complex $\lambda$ as well, with a few adjustments.
Given $\lambda \in \bc$, we will fix $\nsL$ and compute
\[
  M_k^{\nsL} := \max_{n_i, \nsR} \abs{[\nsL, n_i, \nsR]}.
\]
As we will see later, this can be determined as
$M_k^{\nsL} = \max_{n_i \in \mathcal{N}_i} M_k^{(\nsL, n_i)}$ where the finite candidate set $\mathcal{N}_i$ depends on $\nsL = (n_1, \cdots, n_{i-1})$.
Then, maximum over all $\ns$ in $(\bar{\bz} \setminus \{0\})^k$ is obtained as $M_k = M_k^{()} = \max_{n_1 \in \mathcal{N}_1} M_k^{(n_1)} = \max_{n_1 \in \mathcal{N}_1} \max_{n_2 \in \mathcal{N}_2 (n_1)} M_k^{(n_1, n_2)} = \cdots$.
This can be determined similar to the exhaustive algorithm for real $\lambda$.

Recall that for $i \geq 2$, one has
\begin{align*}
  [\ns] & = [n_1, \cdots, n_{i-1}] \cdot \frac{n_i + [n_{i-1}, \cdots, n_2] + [n_{i+1}, \cdots, n_k]}{n_i + [n_{i-1}, \cdots, n_1] + [n_{i+1}, \cdots, n_k]} = [\nsL] \left( 1 - \frac{D}{n_i - x_{\infty}} \right)
\end{align*}
by denoting $D := [n_{i-1}, \cdots, n_1] - [n_{i-1}, \cdots, n_2]$; $x_{\infty}$ is the singular point given as $x_\infty = - ([n_{i-1}, \cdots, n_1] + [n_{i+1}, \cdots, n_k])$.
Remark that the singular point is a complex number, so it might not lie on $\br$.
Then for fixed $\nsL$, one has
\[
  M_k^{\nsL} = \max_{n_i, \nsR} \abs{[\ns]} = \abs{[\nsL] D} \cdot \max_{n_i, \nsR} \left| \frac{1}{D} - \frac{1}{n_i - x_{\infty}} \right|.
\]
Meanwhile, if $i = 1$, we may take $x_{\infty} = - [n_2, \cdots, n_k]$ and $D = \infty$ so that
\[
  M_k^{()}
  = \max_{n_1, \nsR} \abs{[n_1, \nsR]}
  = \max_{n_1, \nsR} \abs{\frac{1 / \lambda}{n_1 + [\nsR]}}
  = \frac{1}{\abs{\lambda}} \max_{n_1, \nsR} \abs{\frac{1}{\infty} - \frac{1}{n_1 - x_{\infty}}}.
\]
In particular, for any $i$, $M_k^{\nsL}$ is obtained when $\abs{D^{-1} - (n_i - x_{\infty})^{-1}}$ is maximal. 

First,
we will compute $x_{\max} \in \bar{\br}$ which maximizes $\abs{[\nsL, x, \nsR]}$.

\begin{lemma}\label{maximal case estimate}
  Let $x_{\infty} = a + b \im$, and suppose $D \: / \: (2 b \im) \neq 0, -1$.
  Then, $(x \mapsto \abs{[\nsL, x, \nsR]}) : \bar{\br} = \br \cup \{ \infty \} \to [0, \infty]$ has a unique local maximum at $x = x_{\max} \in \bar{\br}$, given as
  \begin{equation*}
    x_{\max} - a = \begin{cases}
      - \frac{\abs{D} \abs{D + 2b \im} - \inprod{D}{D + 2b \im}}{2 \Re(D)}, & D \neq \infty, \Re(D) \neq 0 \\
      \infty, & D / (2b \im) \in (-1, 0) \\
      0, & D / (2b \im) \in \bar{\br} \setminus (-1, 0)
    \end{cases},
  \end{equation*} where $\inprod{-}{-}$ is a $\br$-bilinear form given by $\inprod{c + d \im}{c' + d' \im} = cc' + dd'$.
\end{lemma}
\begin{proof}
Consider that the M\"{o}bius transformation $\varphi(z) = z^{-1}$ sends $\bar{\br} - x_{\infty}$ to a line or a circle. When it is mapped to a line, $\bar{\br} - x_{\infty}$ should pass through the origin,
so $x_{\max} = x_{\infty} \in \br$ obtains the maximum $\infty$; one can check this matches the desired conclusion.
Hence, it is enough to show for the general case where $\varphi(\bar{\br} - x_{\infty})$ is a circle passing through the origin.

  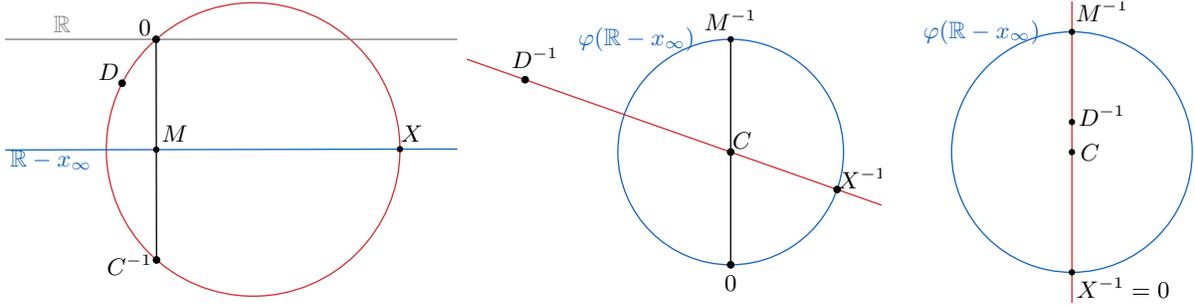
\begin{figure}[H]
    \centering
    \begin{tikzpicture}[line cap=round,line join=round,scale=0.5]
      \clip(-9,-4) rectangle (3,4.5);
      \draw [line width=0.5pt,color=geored] (-2.42,0.07) circle (3.9040107581818986cm);
      \draw [line width=0.5pt] (-5,3)-- (-4.988339894200145,-2.8702262137223373);
      \draw [line width=0.5pt, color=gray] (-10,3) -- (3,3);
      \draw [line width=0.5pt,color=geoblue,domain=-9:3] plot(\x,{(--0.43913329099621135--0.011660105799855103*\x)/5.870226213722337});
      \begin{scriptsize}
        \draw [fill=black] (-5,3) circle (2.5pt);
        \draw[color=black] (-5.3,3.3) node {$0$};
        \draw[color=gray] (-7.5, 3.4) node {$\br$};
        \draw [fill=black] (-4.988339894200145,-2.8702262137223373) circle (2.5pt);
        \draw[color=black] (-5.7,-3) node {$C^{-1}$};
        \draw[color=geoblue] (-7.8,-0.2) node {$\br-x_{\infty}$};
        \draw [fill=black] (-4.9941699471000724,0.06488689313883132) circle (2pt);
        \draw[color=black] (-4.5,0.5) node {$M$};
        \draw [fill=black] (1.4840030566869966,0.07775457146396089) circle (2pt);
        \draw[color=black] (1.8,0.5) node {$X$};
        \draw [fill=black] (-5.903419240287578,1.8326940733985344) circle (2.5pt);
        \draw[color=black] (-6.26,2.14) node {$D$};
      \end{scriptsize}
    \end{tikzpicture}
    \begin{tikzpicture}[line cap=round,line join=round,scale=0.5]
      \clip(1,-3) rectangle (12,5);
      \draw [line width=0.5pt,color=geoblue] (8,1) circle (3cm);
      \draw [line width=0.5pt,color=geored,domain=-4.54:21.06] plot(\x,{(-10.800708883665594--0.9963753582717338*\x)/-2.8297060174917235});
      \draw [line width=0.5pt] (8,4)-- (8,-2);
      \begin{scriptsize}
        \draw [fill=black] (8,1) circle (2.5pt);
        \draw[color=black] (8.3,1.3) node {$C$};
        \draw [fill=black] (8,-2) circle (2.5pt);
        \draw[color=black] (8,-2.5) node {$0$};
        \draw[color=geoblue] (5.5,4) node {$\varphi(\br-{x_{\infty}})$};
        \draw [fill=black] (10.829706017491723,0.0036246417282661536) circle (2.5pt);
        \draw[color=black] (11.5,0.3) node {$X^{-1}$};
        \draw [fill=black] (8,4) circle (2pt);
        \draw[color=black] (8,4.5) node {$M^{-1}$};
        \draw [fill=black] (2.5345287680903645,2.9244617013766323) circle (2.5pt);
        \draw[color=black] (2.8,3.5) node {$D^{-1}$};
      \end{scriptsize}
    \end{tikzpicture}
    \begin{tikzpicture}[line cap=round,line join=round,scale=0.4]
      \clip(-7,-5) rectangle (3,5);
      \draw [line width=0.5pt,color=geoblue] (-1,0) circle (4cm);
      \draw [line width=0.5pt,color=geored] (-1,-6.5) -- (-1,6.5);
      \begin{scriptsize}
        \draw [fill=black] (-1,0) circle (2.5pt);
        \draw[color=black] (-0.4,0) node {$C$};
        \draw [fill=black] (-1,-4) circle (2.5pt);
        \draw[color=black] (0.7,-4.5) node {$X^{-1} = 0$};
        \draw[color=geoblue] (-4,4) node {$\varphi(\br - x_{\infty})$};
        \draw [fill=black] (-1,4) circle (2.5pt);
        \draw[color=black] (0,4.7) node {$M^{-1}$};
        \draw [fill=black] (-1,1) circle (2.5pt);
        \draw[color=black] (0,1.2) node {$D^{-1}$};
      \end{scriptsize}
    \end{tikzpicture}

    \caption{Left: $\br - x_{\infty}$ and relevant points before transformation. \\
    Middle, Right: After transformation $\varphi$;
    figure on right depicts the special case where $D^{-1}$ lies on $\im \br$ above $C$.}
    \label{fig:circles}
  \end{figure}

  Note we can identify the circle's center $C \in \bc$ as follows.
  Let $M = - \Im(x_{\infty}) \im$, an intersection of the imaginary line $\im \br$ and $\br - x_{\infty}$.
  $\varphi$ sends $M$ to an intersection of the imaginary line and $\varphi(\br - x_{\infty})$,
  where the two meets orthogonally.
  Hence, $C$ is the midpoint of $M^{-1}$ and $0$, and so $C^{-1} = 2M = - 2 \Im(x_{\infty}) \im$.

  Now, observe that the value to maximize, $\abs{D^{-1} - (x - x_{\infty})^{-1}}$, is the distance from $(x - x_{\infty})^{-1} \in \varphi(\br - x_{\infty})$ to $D^{-1}$.
  Clearly it has a unique local maximum,
  which is obtained when $D^{-1}$, the center $C$, and $(x - x_{\infty})^{-1}$ are collinear in this order.
  Remark that the same logic works for $D = \infty$ as well, since $D^{-1} = 0 \neq C$.
  Let us denote the local maximum as $x_{\max}$, and let $X = x_{\max} - x_{\infty}$.
  Then, one can calculate $X-M = x_{\max} - a$ from the fact that $X^{-1}, C, D^{-1}$ are collinear,
  using either cross ratio or basic trigonometry.
  This gives the desired result, where the bilinear form originates from $\Re(\bar{z}w) = \inprod{z}{w}$.
\end{proof}

\begin{remark}
  $\abs{[\nsL, x, \nsR]}$ is constant if $D = 0$ or $D = - 2 \Im(x_{\infty}) \im$.
  Indeed, one has $[\nsL, x, \nsR] = [\nsL]$ if $D = 0$,
  and the center $C$ coincides with $D^{-1}$ if $D = -2 \Im(x_{\infty}) \im$.
\end{remark}

Now, we can show the statement analogous to Proposition \ref{max wcf}.

\begin{proposition}\label{maximal ni estimate}
  The maximum $M_k^{\nsL}$ can be obtained by $\ns = (\nsL, n_i, \nsR)$ satisfying
  \[
    \begin{cases}
      \abs{n_i - \left(\Re(x_{\infty}) - \frac{\delta(2 \Im(x_{\infty}))}{2 \Re(D)} \right)} < 2, & \Re(D) \neq 0, \\
      \abs{n_i - \Re(x_{\infty})} < 2 \text{ or } n_i = \infty, & \Re(D) = 0 \text{ or } D = \infty
    \end{cases}
  \]
  where $\delta(t) = \abs{D} \abs{D + t \im} - \inprod{D}{D + t \im}$.

\end{proposition}
\begin{proof}
 This can be checked in case-by-case manner. The inequality comes from $\abs{n_i - x_{\max}} < 2$ in cases where $n_i$ is finite, whose reasoning is analogous to the proof of Proposition \ref{max wcf}.
\end{proof}

\begin{remark}
  The function $\delta(t) = \abs{D} \abs{D + t \im} - \inprod{D}{D + t \im}$ introduced in Proposition \ref{maximal ni estimate} only has a local minimum of $0$ at $t = 0$, without any local maximum; this can be checked by looking at the derivative.
  Hence, $\delta(t)$ is maximized at the boundary of closed interval.
\end{remark}

Now, we are ready to compute the maximum $M^{\nsL} = \max_{n_i, \nsR} \abs{[\nsL, n_i, \nsR]}$ given $\nsL$,
whose process is outlined below.
\begin{enumerate}
  \item Compute $D$, and compute range of $x_{\infty}$ (denoted $\mathcal{B} \subset \bc$)
        using $\abs{x_{\infty} + [\overleftarrow{\nsL}]} \leq \max \abs{[\nsR]} = M_{k-i}$.
  \item Compute range $\mathcal{N}_i \subset \bar{\bz} \setminus \{0\}$ of $n_i$ depending on $D$.
  \begin{enumerate}
    \item If $\Re(D) = 0$, $\abs{n_i - \Re(x_{\infty})} < 2$ or $n_i = \infty$,
          so $\mathcal{N}_i = (\Re(\mathcal{B}) + (-2, 2)) \cup \{\infty\}$.
    \item If $\Re(D) \neq 0$, compute $\Delta = \max(\delta(a), \delta(b))$ for $2\Im(\mathcal{B}) = [a, b]$.
          By the remark above, one has $\delta(2 \Im(x_{\infty})) \in [0, \Delta]$. Together with Proposition \ref{maximal ni estimate}, one has
          \[ \mathcal{N}_i = \Re(\mathcal{B}) - \frac{1}{2 \Re(D)} \cdot [0, \Delta] + (-2, 2). \]
  \end{enumerate}
  \item Compute $M_k^{\nsL} = \max_{n_i \in \mathcal{N}_i} M_k^{(\nsL, n_i)}$ recursively.
        Remark that for $n_i = \infty$, one always have
        $M^{(\nsL, \infty)} = \abs{[\nsL]}$ as $[\nsL, \infty, \nsR] = [\nsL]$.
        Thus, we do not have to further compute on the case $n_i = \infty$.
\end{enumerate}

The desired maximum $M_k = M_k^{()}$ can be obtained by the above process in a Depth-First Search(DFS) fashion, similar as the exhaustive algorithm for real $\lambda$; we begin by choosing $n_1, n_2, \cdots$ from candidates $n_i \in \mathcal{N}_i (n_1, \cdots, n_{i-1})$, descending down to $n_k$; then, each maximum $M_k^{\nsL}$ will be computed from the computed $M_k^{\nsL, n_i}$'s as we recurse back up.
The $M_k$ obtained as above can be used to determine the minimal $u$-degree.

\begin{remark}
   Unlike the case when $\lambda$ is limited to the real numbers, $M_k$ might not be strictly increasing, but it is still increasing. Specifically, a maximum can be obtained as $M_k = [\nsL, \infty, \nsR]$.
   For instance, one has $M_1 = \wcf{1/\lambda}{1, \infty} = M_2$ for $\lambda = \sqrt{-1}$ case.
\end{remark}

Note that the prescribed algorithm is generally slower than the algorithm \ref{max wcf alg}. This is because the range of $x_{\infty}$, $\mathcal{B}$, is way larger on complex plane than on real line. As a result, the candidate set $\mathcal{N}_i$ of $n_i$ is larger in this case.
Nevertheless, this algorithm was able to determine the following relation numbers in a few seconds.

\begin{itemize}
    \item $\lambda = \im$ has minimal $u$-degree 5; $\ns = (1, 1, -1, -1, 1)$
    \item $\lambda = \frac{1 + \im}{2}$ has minimal $u$-degree 4; $\ns = (1, -1, -2, 1)$
    \item $\lambda = \frac{1 + \sqrt{-3}}{2}$ has minimal $u$-degree 4; $\ns = (-1, 1, 1, -1)$
    \item $\lambda = \frac{- 1 + \sqrt{-3}}{2}$ has minimal $u$-degree 4; $\ns = (-1, -1, 1, 1)$
\end{itemize}

We have also devised an improved version of this algorithm by giving separate bounds to real parts and imaginary parts of $[\ns]$, whose details will be omitted. The following relation numbers were found in this way.
\begin{itemize}
    \item $\lambda = \frac{2 + \im}{3}$ has minimal $u$-degree 5; $\ns = (-3, 1, -3, -1, 1)$
    \item $\lambda = \frac{3 + 4 \im}{5}$ has minimal $u$-degree 5; $\ns = (-5, 1, -1, -1, 1)$
    \item $\lambda = \frac{2 + \sqrt{-5}}{3}$ has minimal $u$-degree 5; $\ns = (-3, 1, -1, -1, 1)$
    \item $\lambda = \frac{1 + \sqrt{-8}}{3}$ has minimal $u$-degree 5; $\ns = (-3, -1, 1, 1, -1)$
\end{itemize}


\nocite{}

\end{document}